\documentclass[11pt]{article}
\usepackage{amssymb}
\usepackage{graphicx}
\usepackage{latexsym}
\usepackage{amsmath,amscd}
\def\part#1{\frac{\partial\phantom{q}}{\partial#1}}

\newenvironment{rmk}{\begin{trivlist}\item[]{\bf Remark:} }
{\end{trivlist}}

\newenvironment{ex}{\begin{trivlist}\item[]{\bf Example:} }
{\end{trivlist}}

\newenvironment{prf}{\begin{trivlist}\item[]{\bf Proof:} }
{\hfill $\Box$ \end{trivlist}}
\newenvironment{lemprf}{\begin{trivlist}\item[]{\bf Proof:} }
 {\end{trivlist}}
\newtheorem{thm}{Theorem}

\newtheorem{prp}[thm]{Proposition}
\newtheorem{lemma}[thm]{Lemma}

\newcommand{\lie}[1]{\mathfrak{#1}}
\def\End{\mathop{\rm End}\nolimits}

\def\ker{\mathop{\rm ker}\nolimits}
\def\coker{\mathop{\rm coker}\nolimits}

\def\grad{\mathop{\rm grad}\nolimits}

\def\tr{\mathop{\rm tr}\nolimits}

\def\Real{\mathop{\rm Re}\nolimits}

\newcommand{\R}{\mathbf{R}}
\newcommand{\C}{\mathbf{C}}
\newcommand{\K}{\mathbf{H}}
\newcommand{\Z}{\mathbf{Z}}

\newcommand{\CP}{{\mathbf P}}

\newcommand{\RP}{{\mathbf R}{\rm P}}
\newcommand{\KP}{{\mathbf H}{\rm P}}
\begin{document}
\title{On the hyperk\"ahler/quaternion K\"ahler correspondence }
 \author{Nigel Hitchin\\[5pt]}
 \maketitle
\centerline{{\it Subject classification}: {Primary 53C26, 53C28} }
\tableofcontents
\section{Introduction}
This paper deals with a correspondence between a hyperk\"ahler manifold $M$ with a circle action which preserves just one complex structure $I$ and a quaternionic K\"ahler manifold $\hat M$   of the same dimension with a circle action. Its origins lie in the physicist's c-map construction but a mathematical theorem in the language of differential geometry was given by A.Haydys \cite{Hay}. Our aim here is to describe this from the viewpoint of twistor theory. The twistor space of a hyperk\"ahler manifold is a holomorphic fibration of complex symplectic manifolds over the projective line, and the twistor space of a quaternionic K\"ahler manifold is a holomorphic contact manifold. Several natural constructions in this area, such as quotients, are more transparent when seen in twistor terms and this is also the case here.

The correspondence is in some ways not symmetrical. To pass from the quaternionic K\"ahler side to the hyperk\"ahler side, one may simply take the hyperk\"ahler quotient of the Swann bundle by the lifted circle action. However, to go in the other direction involves the introduction of a natural hyperholomorphic line bundle on the  hyperk\"ahler manifold. Most of the paper will actually be concerned with this new feature, one which is particularly  interesting when considering the hyperk\"ahler moduli spaces of  gauge-theoretic equations. Given the line bundle, the correspondence consists of lifting the geometric circle  action on the hyperk\"ahler manifold $M$ to its principal $U(1)$-bundle and taking $\hat M$ to be the quotient manifold. Then $U(1)$ induces a geometrical action on $\hat M$. 

To describe a hyperk\"ahler metric in twistorial terms requires two things --  the twistor space $Z$ itself, and a space of rational curves, sections of the fibration, called the twistor lines. In many cases we know the twistor space but have no analytic expression for the hyperk\"ahler metric or the twistor lines. However, a hyperholomorphic line bundle on $M$ is uniquely determined by a holomorphic line bundle $L_Z$ on the twistor space. One of the positive features about our description of it is that $L_Z$ only depends on the geometry of the twistor space and not on the twistor lines.  Thus, in examples, we can sometimes describe it by a single holomorphic  transition function even when the metric is not explicit. 

Haydys introduced the line bundle via its curvature form $\omega_1+dd^c_1\mu$, where $\omega_1$ is the invariant K\"ahler form and $\mu$ the moment map. In this article we first define a natural holomorphic line bundle  on the twistor space of a hyperk\"ahler manifold with a circle action, and then show that it agrees with the differential geometric construction of Haydys. The holomorphic description uses a \v Cech cohomological approach but one may interpret this by saying that $L_Z$ on $Z\rightarrow \CP^1$ has a meromorphic connection with poles on the fibres over $0$ and $\infty$. The  curvature of this connection is a closed meromorphic 2-form which restricts on each fibre over $\CP^1\backslash\{0,\infty\}$ to a multiple of the holomorphic symplectic form.  One way to view this is to say that the twistor space, seen as a fibration over the projective line, is a symplectic manifold defined over the field of functions in one variable and the line bundle is then the  quantum line bundle over this field.

 Following the twistor construction we look at various examples including Taub-NUT space, the semi-flat case, the Legendre transform construction of \cite{HKLR}, the cotangent bundle construction of \cite{Feix1} and monopole moduli spaces. The line bundle can in these cases be economically defined by transition functions.  The Higgs bundle moduli spaces require a different description: ours is in terms of  determinant lines and zeta-function determinants. 
 
 Finally we describe how to pass back and forth between the twistor spaces of the  corresponding hyperk\"ahler and quaternionic K\"ahler manifolds. The contact form for the quaternionic K\"ahler metric is defined in terms of the connection form for the meromorphic connection on $L_Z$. 
  
The author wishes to thank Andrew Swann for introducing him to the subject  at the Aarhus  QGM Centre, and Simon Salamon, Martin Ro\v cek and Stefan Vandoren for useful discussions. He also thanks the Simons Center for Geometry and Physics for support while this research was being carried out. 

\section{The hyperholomorphic line bundle}
\subsection{The differential geometric approach}\label{diff}
This is the approach of Haydys \cite{Hay}. Let $M$ be a hyperk\"ahler manifold of dimension $4k$ with K\"ahler forms $\omega_1,\omega_2,\omega_3$ corresponding to complex structures $I,J,K$. Suppose there is a circle action which fixes $\omega_1$ and rotates $\omega_2$ and $\omega_3$. 

\begin{ex}\label{ex1}
A model for this is  $\C^2$ with the flat metric and 
 $$\omega_1=\frac{i}{2}(dz \wedge d\bar z+dw \wedge d\bar w), \qquad \omega_2+i\omega_3=dz\wedge dw.$$
Then  $(z,w)\mapsto (z,e^{i\theta}w)$ takes  $\omega_2+i\omega_3$ to $e^{i\theta}(\omega_2+i\omega_3).$
\end{ex}
The action  generates a vector field $X$ such that 
\begin{equation}
{\mathcal L}_X\omega_1=0,\quad {\mathcal L}_X\omega_2=-\omega_3,\quad {\mathcal L}_X\omega_3=\omega_2.
\label{Xomega123}
\end{equation}
Assume there is a moment map for the action on $\omega_1$ then $i_X\omega_1=d\mu$  and, as shown in \cite{HKLR}, for any vector field $Y$
$$Kd\mu(Y)=\omega_1(X,KY)=g(IX,KY)=-g(KIX,Y)=-g(JX,Y)=-i_X\omega_2(Y),$$
from which
$$dKd\mu=-d(i_X\omega_2)=-{\mathcal L}_X\omega_2=\omega_3.$$
For a complex structure $I$, the real operator $d^c$ is defined as $I^{-1}dI$, so on $M$ we have three such operators $d^c_1,d^c_2,d^c_3$ and hence we can write the above as 
$dd^c_3\mu=-\omega_3$ and  $dd^c_2\mu=-\omega_2$.

We now have from \cite{Hay},

\begin{prp} $\omega_1+dd^c_1\mu$ is of type $(1,1)$ with respect to  $I,J$ and $K$.
\end{prp}

\begin{prf} First recall that because the Levi-Civita connection is torsion-free ${\mathcal L}_XY=\nabla_XY-(\nabla X)(Y)$ where $\nabla X$ is a section of $T\otimes T^*=\End T$. The Lie derivative  on a differential form $\alpha$ is similar:  ${\mathcal L}_X\alpha=\nabla_X\alpha-(\nabla X)(\alpha)$ where $\nabla X$ acts through the action of $\End T$ on forms. Since $M$ is hyperk\"ahler, $\nabla_X$ acts trivially on $\omega_1,\omega_2,\omega_3$  and by hypothesis ${\mathcal L}_X$ only acts trivially on $\omega_1$ but preserves the space spanned by these forms. We deduce that at each point $\nabla X$ lies in the subbundle of $\End T$ corresponding to the Lie algebra  $\lie{u}(1)+\lie{sp}(n)$.

Using the metric to identify $2$-forms with skew-adjoint endomorphisms of $T$, the subspace $\lie{sp}(n)$ is the space of $2$-forms of type $(1,1)$ with respect to $I,J,K$ and the trivial $\lie{u}(1)$ factor consists of multiples of the K\"ahler form $\omega_1$. 
From  $IX=\grad \mu$ we see that  $\nabla X\in \End T$ is equivalent to  $dd^c_1\mu\in \Lambda^{1,1}T^*$, so to prove the Proposition we must show that the $\lie{u}(1)$ component of $dd^c_1\mu$ is the same as that of $-\omega_1$. If we introduce the Lefschetz operators $\Lambda_1,\Lambda_2,\Lambda_3: \Lambda^pT^*\rightarrow \Lambda^{p-2}T^*$ this is equivalent to proving $\Lambda_1\omega_1=-\Lambda_1dd_1^c\mu$.

Now $\Lambda_1dd_1^c\mu=-\Delta \mu$ where $\Delta$ is the Laplacian and 
$$2k=\Lambda_1\omega_1=\Lambda_2\omega_2=-\Lambda_2dd^c_2\mu =\Delta\mu$$
using  $dd^c_2\mu=-\omega_2$. Hence $\Lambda_1\omega_1=-\Lambda_1dd_1^c\mu$ as required.
\end{prf}

\begin{ex} Consider $\C^2$  as above. The circle action $(z,w)\mapsto (z,e^{i\theta}w)$ gives the vector field 
 $$X=iw\frac{\partial}{\partial w}-i\bar w\frac{\partial}{\partial \bar w}$$
 and the moment map 
 $$\mu=-\frac{1}{2}\vert w\vert^2 \qquad d\mu=-\frac{1}{2} (w d\bar w+\bar w dw).$$
So $dd_1^c\mu=-idw\wedge d\bar w$ and  $$\omega_1+dd_1^c\mu=\frac{i}{2}(dz\wedge d\bar z-dw \wedge d\bar w).$$
 Note that in  four dimensions the condition on a form to be of type $(1,1)$ with respect to  $I,J$ and $K$ is   equivalent to it being anti-self-dual.
 \end{ex}

If the cohomology class of $\omega_1/2\pi$ is integral (and hence also that of $F=\omega_1+dd^c_1\mu$) then we may regard $F$ as the curvature of a $U(1)$-connection on a principal circle bundle $P$. If $M$ is simply connected this is unique up to gauge equivalence. (Note that this is the quantum line bundle of the symplectic form $\omega_1$ but with Hermitian metric rescaled by $e^{\mu}$). We need to lift the circle action to $P$. A lift of the vector field $X$ is given by an invariant  function $\phi$ such that $i_XF=d\phi$. 

In our case
$$i_XF=d\mu+i_Xdd_1^c\mu=d\mu+{\mathcal L}_Xd_1^c\mu-d(i_Xd_1^c\mu)$$
using the Cartan formula for Lie derivative of a form. But $X$ preserves the complex structure $I$ and the moment map $\mu$ so ${\mathcal L}_Xd_1^c\mu=d_1^c({\mathcal L}_X \mu)=0$, which means that
$$i_XF=d(\mu+i_XId\mu).$$ 
However $i_XId\mu=d\mu(IX)=g(IX,IX)=g(X,X)$ so, given $\mu$, we have a natural choice of lift given by $\phi=\mu+g(X,X)$.

To lift the circle action we require that the equivariant cohomology class $[\omega+u\phi]$ should be integral. Any two liftings differ by a homomorphism from the circle to the $U(1)$ acting on the fibres of the principal bundle. This involves changing the function $\phi$, or equivalently $\mu$,  by an integer.

A connection whose curvature is of type $(1,1)$ with respect to $I,J,K$ is called {\it hyperholomorphic}. Since the $(1,1)$ forms for $I$ are precisely those $2$-forms on which the Lie algebra action of $I $ is trivial, to be $(1,1)$ with respect to $I,J,K$ implies $(1,1)$ with respect to a complex structure which is a linear combination of $I,J,K$ -- any of the $2$-sphere family of complex structures on a hyperk\"ahler manifold. Given the integrality constraint, we have just produced a hyperholomorphic line bundle on $M$. 

\begin{rmk} There is one standard method of producing hyperholomorphic bundles and this is where the hyperk\"ahler manifold is expressed as a hyperk\"ahler quotient. Recall that if $G$ acts on $M$ preserving $\omega_1,\omega_2,\omega_3$ and  we have a hyperk\"ahler moment map $\mu:M\rightarrow {\lie g}^*\otimes \R^3$, then $\bar M=\mu^{-1}(0)/G$ has an induced hyperk\"ahler metric \cite{HKLR}. The space  $\mu^{-1}(0)\subset M$ is a principal $G$-bundle over $\bar M$ and the induced metric from $M$ defines horizontal subspaces to the tangent spaces  of the $G$-orbits. This defines a connection which is hyperholomorphic.  
\end{rmk}
\subsection{The twistor approach} \label{twist}
To each hyperk\"ahler manifold $M$ of dimension $4k$ one may associate its twistor space, a complex manifold $Z$ of (complex) dimension $2k+1$. Differentiably it is the product $M\times S^2$ but with complex structure at $(m,{\bf x})$ defined by $(I_{\bf x},I)$ where $I$ is the standard complex structure on $S^2\cong \CP^1$ and $I_{\bf x}=x_1I+x_2J+x_3K, {\bf x}=(x_1,x_2,x_3)\in S^2$. A hyperholomorphic connection has curvature of type $(1,1)$ with respect to all $I_{\bf x}$ and it follows that the connection on $M$ pulled back to the product $M\times S^2$ defines a holomorphic structure on the bundle. To interpret the hyperholomorphic line bundle of the previous section in twistorial terms we are therefore seeking a holomorphic line bundle $L_Z$ on its twistor space. If it is trivial on each twistor line this is also sufficient, by the hyperk\"ahler version of the Atiyah-Ward correspondence.

As a complex manifold, the twistor space has the following features \cite{HKLR}:
\begin{itemize}
\item
a holomorphic projection $\pi:Z\rightarrow \CP^1$
\item
a holomorphic section $\omega$ of $\Lambda^2T^*_F(2)$ (where $T_F$ is the tangent bundle along the fibres and the $(2)$ denotes the tensor product with $\pi^*{\mathcal O}(2)$) such that $\omega$ defines a holomorphic symplectic form in each fibre
\item
a real structure preserving this data and inducing the antipodal map on $\CP^1$.
\end{itemize}

We shall produce, using \v Cech cocycles,  a holomorphic line bundle on $Z$ from the data of the circle action and then in the next section show that it coincides with the bundle defined by the procedure of Section \ref{diff}.

It is useful to suppress for the moment the integrality condition on $\omega_1/2\pi$.  There is still  a geometric way to phrase this, which is useful in the holomorphic context.   If $P$ is a holomorphic principal $\C^*$-bundle over a complex manifold $X$ then $TP/\C^*$ is a bundle on $X$ whose local sections are the invariant vector fields on $P$.  More generally a bundle $E$ which is an extension
$$0\rightarrow {\mathcal O}\rightarrow E\rightarrow T\rightarrow 0$$
with a Lie bracket operation compatible with the projection to $T$ is called a {\it Lie algebroid}. With respect to local trivializations it is defined by a $1$-cocycle of closed $1$-forms.  The exact cohomology sequence of the sequence of sheaves $0\rightarrow \C\rightarrow {\mathcal O}\rightarrow d{\mathcal O}\rightarrow 0$ associates  to each such $1$-cocycle a class in $H^2(X,\C)$ which is the  characteristic class of the algebroid. If this lies in $H^2(X,2\pi i\Z)$ then we can replace this with a   cocycle of the form  $g_{UV}^{-1}dg_{UV}$ for the transition functions $g_{UV}$ of a principal  $\C^*$-bundle and then the Lie algebroid is $TP/\C^*$ (usually called an Atiyah algebroid).  Although we are looking for a holomorphic line bundle on the twistor space $Z$, it is more convenient to seek a holomorphic Lie algebroid first and then discuss the integrality condition later. 

\begin{prp} \label {alg} On each fibre of the twistor space $Z$ of a hyperk\"ahler manifold there exists a natural holomorphic Lie algebroid.
\end{prp}

\begin{prf} The twistor space $\pi: Z\rightarrow \CP^1$ is covered by open sets $\{U, V,\dots\}$ where each set is the product of an open set in $\C^{2k}$ and an open set in $ \CP^1$, with affine parameter $\zeta$. With respect to the product structure,  the local vector field $d/d\zeta$ on $\CP^1$ lifts to a unique vector field $Z_U$ on $U$. On the intersection of two such open sets $U,V$ we have $Z_V-Z_U=X_{UV}$, which is a vector field along the fibres of $U\cap V\subset Z\rightarrow \CP^1$. Invariantly speaking $d/d\zeta$ is a local trivialization of the tangent bundle of $\CP^1$ which is isomorphic to ${\mathcal O}(2)$ and then $X_{UV}$ is a $1$-cocycle with values in $T_F(-2)$.  It  defines the cotangent bundle of $Z$ as an extension:
$$0\rightarrow {\mathcal O(-2)} \rightarrow T_Z^*\rightarrow T_F^*\rightarrow 0.$$
Now take $\omega$, which is a section $\Lambda^2T^*_F(2)$, and form the cocycle $\theta_{UV}=i_{X_{UV}}\omega\in T^*_F$. This defines an extension 
$$0\rightarrow {\mathcal O} \rightarrow E\rightarrow T_F\rightarrow 0.$$
This extension is defined by a cocycle of 1-forms but for a Lie algebroid we need closed forms, which means choosing $Z_U$ more carefully. This local vector field projects to one on $\CP^1$ and so when integrated takes fibres to fibres. It also acts on the tangent bundle ${\mathcal O}(2)$ and so we have a Lie derivative action on $\Lambda^2T^*_F(2)$. We ask that ${\mathcal L}_{Z_U}\omega=0$. This is possible because by Darboux' s theorem locally $Z$ is a product of a symplectic manifold and an open set in $\C$. With this choice we see that ${\mathcal L}_{X_{UV}}\omega=0$, so on a fixed fibre $di_{X_{UV}}\omega= d\theta_{UV}=0$.This cocycle defines the Lie algebroid structure on $E$ for each fibre $\zeta=const$. 
\end{prf}
 \vskip .25cm
The existence of the hyperholomorphic line bundle requires also an $S^1$-action on $M$. This extends naturally to $Z$ preserving the holomorphic structure and therefore defines a holomorphic vector field $Y$ on $Z$. The original action induced a rotation on the 2-sphere of complex structures, leaving fixed two points which we may assume are  $\zeta=0$ and $\zeta =\infty$. Then the vector field $Y$ projects to the vector field $i\zeta d/d\zeta$ 
on $\CP^1$, and the section $\omega$ of $\Lambda^2T^*_F(2)$ is invariant.  Invariance again means that  $Y$ preserves the fibration and the line bundle ${\mathcal O}(2)$ on the quotient, and so there is a natural Lie derivative.
Thus, regarding $Z_U$ as a section of $T_Z(-2)$, we have ${\mathcal L}_YZ_U=W_U$ where $W_U$ is a local section of $T_F(-2)$ and 
${\mathcal L}_YX_{UV}=W_V-W_U.$ Note that if $U$ is invariant under the circle action then by averaging over the circle we can choose $Z_U$ to be invariant and then $W_U=0$.

The vector fields $X_{UV}$ and $W_U$ preserve the symplectic forms along each fibre hence 
\begin{equation}
{\mathcal L}_Y(i_{X_{UV}}\omega)=i_{W_V}\omega-i_{W_U}\omega=d_Ff_V-d_Ff_U
\label{LUV}
 \end{equation}
for locally defined functions $f_U$.  

Consider now the fibre $Z_0$ given by $\zeta=0$: this is the manifold $M$ with complex structure $I$. As noted above, since $d_F(i_{X_{UV}}\omega)=0$, $i_{X_{UV}}\omega$  is closed on each fibre and in particular $Z_0$. But  $Y$ is tangential to $Z_0$ and   it follows from (\ref{LUV}) therefore that 
$$d(i_Yi_{X_{UV}}\omega-f_V+f_U)=0.$$
 This defines a $1$-cocycle with values in $\C$ and so if $H^1(M,\C)=0$ (which we assume for the remainder of the paper) then we may choose the $\{f_U\}$ so that  
 \begin{equation}
 f_V=f_U+i_Yi_{X_{UV}}\omega.
 \label{fdefine}
 \end{equation}
 
 Using the isomorphism $T_Z^*\cong {\mathcal O}(-2)\oplus T^*_F$ provided by the local product decomposition on $U$ a form $\varphi$ along the fibres gives a form $\varphi_U$ on $U$  characterized by  $i_{Z_U}\varphi_U=0$. Applying this to $\omega$  we have sections $\omega_U$ of $\Lambda^2T_Z^*(2)$ where 
$ \omega_V-\omega_U=i_{X_{UV}}\omega\wedge d\zeta$.
 It follows from (\ref{fdefine})  that on $U\cap V$ 
$$(i_Y\omega_V-f_Vd\zeta)-(i_Y\omega_U-f_Ud\zeta)=i_Yi_{X_{UV}}\omega d\zeta-i_Yi_{X_{UV}}\omega d\zeta=0.$$
which means that $i_Y\omega_U-f_Ud\zeta$ is the restriction  of  a global section $\varphi_0$ of $T_Z^*(2)$
on $Z_0$ to $U\cap Z_0$.

\begin{rmk}   Note that $\varphi_0$ is not quite unique, for we can add  the same constant to each $f_U$, but the restriction of $\varphi_0$ to $T^*_F\cong T^*_M$ is the canonical 1-form $i_Y\omega$.   Using the local splittings in the computations above, we observe that $\tilde Y\vert_U=-f_U +Y$ is a lift of the vector field $Y$ to the algebroid $E$, and changing $f_U$ by a constant gives a different lift.  
\end{rmk}

Now apply the real structure to transform $\varphi_0$ to a form $\varphi_{\infty}$ over the fibre $Z_{\infty}$. The divisor  $D=Z_0+Z_{\infty}$ is the zero set  of the section $s=p^*(i\zeta d/d\zeta)$ of ${\mathcal O}(2)$ and so in the cohomology sequence of the exact sequence of sheaves
$$0\rightarrow T_Z^*\stackrel {s}\rightarrow T_Z^*(2)\rightarrow T_Z^*(2)\vert_D\rightarrow 0$$
$\varphi=\varphi_0+\varphi_{\infty}\in H^0(D, T_Z^*(2))$ maps to an element  $\delta(\varphi) \in H^1(Z,T_Z^*)$ and defines on $Z$ an extension
$$0\rightarrow {\mathcal O} \rightarrow E\rightarrow T_Z\rightarrow 0.$$
\begin{prp} \label{algebroid}The extension above has the structure of a Lie algebroid which restricts on each fibre to the algebroid of Proposition \ref{alg}.
\end{prp}
\begin{prf}
We shall find a cocycle of closed forms to define the extension.  First recall the definition  of the   homomorphism $\delta:H^0(D,T_Z^*(2))\rightarrow H^1(Z,T_Z^*)$ in the exact cohomology sequence: on an open set $U$ on $Z$, extend $\varphi\vert_U$ on $D\cap U$ to $\varphi_U$ on $U$ and then the cocycle $(\varphi_V-\varphi_U)/s$ in $T_Z^*$ represents the class $\delta(\varphi)$. Using the local product decomposition on $U$ we can write $Y=Y_U+i\zeta Z_U$ where $Y_U$ is a local vector field along the fibres. Then
$$\varphi_U=i_{Y_U}\omega_U-f_Ud\zeta$$ is a local extension. So the cocycle is 
$(\varphi_V-\varphi_U)/s$ 
and we need to show this is closed.

From the definition of $f_U$ we have 
$$d(f_Ud\zeta/s)=i_{W_U}\omega/s\wedge d\zeta=i_{W_U}\omega_U/s\wedge d\zeta=i_{[Y,Z_U]}\omega_U/s\wedge d\zeta=-i_{Z_U}{\mathcal L}_Y(\omega_U/s)\wedge d\zeta$$
since $i_{Z_U}\omega_U=0$. Now $\omega/s$ is invariant  by the action so 
${\mathcal L}_Y(\omega/s\wedge d\zeta)=(\omega/s)\wedge {\mathcal L}_Yd\zeta=i(\omega/s)\wedge d\zeta$.  But $\omega_U\wedge d\zeta=\omega\wedge d\zeta$ which means that  ${\mathcal L}_Y(\omega_U/s)=\beta\wedge d\zeta$ for some $\beta$ which we can take to satisfy $i_{Z_U}\beta=0$.
Thus 
$$i_{Z_U}{\mathcal L}_Y(\omega_U/s)\wedge d\zeta=-\beta\wedge d\zeta=-{\mathcal L}_Y(\omega_U/s)$$
and  so $d(f_Ud\zeta/s)={\mathcal L}_Y(\omega_U/s)$. Then using the Cartan formula 
\begin{equation}
d(\varphi_U/s)=d(i_{Y_U}\omega_U/s)-{\mathcal L}_Y(\omega_U/s)=-i_Yd(\omega_U/s).
\label{Fequation}
\end{equation}
From this we have 
$$d(\varphi_V-\varphi_U)/s=i_Yd((\omega_U-\omega_V)/s)=i_Yd(i_{X_{UV}}\omega/s\wedge d\zeta).$$
But $i_{X_{UV}}\omega$ is closed on each fibre and hence $d(i_{X_{UV}}\omega/s\wedge d\zeta)=0$ and the cocycle is closed. 

 Writing
$i_{Y_V}\omega_V-i_{Y_U}\omega_U=i_{Y_U}(\omega_V-\omega_U)+i_{(Y_V-Y_U)}\omega_U$ we obtain 
$$i_{Y_V}\omega_V-i_{Y_U}\omega_U=i_{Y_U}\theta_{UV}d\zeta +i\zeta i_{X_{UV}}\omega_U.$$ 
Restrict the cocycle  to a fibre $\zeta=c\ne 0$ and we therefore obtain
$$(\varphi_V-\varphi_U)/s=i_{X_{UV}}\omega$$
 the cocycle for the Lie algebroid of Proposition \ref{alg}. By continuity this also holds for $\zeta=0,\infty$.
\end{prf}

  From this we obtain the following:

\begin{thm} \label{twistline} Let $Z$ be the twistor space of a hyperk\"ahler manifold $M$  with a circle action as above. Suppose that the canonical Lie algebroid of Proposition 2 on the fibre at $\zeta=0$ is defined by a principal $\C^*$-bundle.  Then there exists a  section  $\varphi$ of $T^*_Z(2)$ on the divisor $D\subset Z$ such that  the Lie algebroid of  Proposition \ref{algebroid}  uniquely defines a line bundle $L_Z$ on $Z$ which is trivial on every twistor line, and which corresponds to  a hyperholomorphic line bundle $L$ on $M$. 
\end{thm}

\begin{prf} The Lie algebroid of Proposition \ref{algebroid}  is defined by  $\delta(\varphi)\in H^1(Z,T^*_Z)$, which  lies in the image of $H^1(Z,d{\mathcal O})$ in $H^1(Z,T^*_Z)$. For uniqueness we need this to be injective. So consider  the exact sequence of sheaves $$0\rightarrow d{\mathcal O}\rightarrow \Omega^1\rightarrow d\Omega^1\rightarrow 0.$$ Restricted to a twistor line in $Z$, the cotangent bundle is $T^*_Z\cong {\mathcal O}(-2)\oplus \C^{2k}(-1)$  so any holomorphic form  vanishes on the line and hence everywhere since there is such a line through each point. In particular  $H^0(Z,d\Omega^1)=0$ and so from the exact cohomology sequence   $H^1(Z,d{\mathcal O})$ injects into  $H^1(Z,\Omega^1)=H^1(Z,T^*_Z)$.

As far as the characteristic class is concerned,   we have the $C^{\infty}$ product $Z=M\times \CP^1$ and so  $H^2(Z,\C)\cong H^2(M,\C)\oplus H^2(\CP^1,\C)$. The second factor is determined by restriction (as a cocycle of $1$-forms) to a twistor line, so this is $\delta\{f_U(0),\bar f_V(\infty)\}\in H^1(\CP^1,{\mathcal O}(-2))\cong \C$. By changing $f_U$ by a constant (and $f_V$ by its conjugate) we can make this class zero.    
The Lie algebroid of Proposition \ref{algebroid} restricts to the canonical algebroid on the fibre so by assumption the integrality condition holds  for the Lie algebroid on $Z$. 

The unitarity of the connection follows from the reality condition for the section of $T^*_Z(2)$ on the real divisor $D$. 
\end{prf}

\vskip .25cm
Our construction reveals a new aspect of the line bundle $L_Z$. Recall that the Atiyah class in $H^1(Z,T^*_Z)$ of a holomorphic line bundle is the obstruction to the existence of a holomorphic connection. In our construction it is of the form $\delta(\varphi)$ for $\varphi$ a section of $T^*_Z(2)$ on the divisor $D$ and as a consequence we have:

\begin{prp} \label{mero} The line bundle $L_Z$ on the twistor space $Z$ admits a meromorphic connection with a simple pole on the divisor $D=\pi^{-1}(0,\infty)$ . It has the following properties:
\begin{enumerate}
\item
its curvature is a closed meromorphic 2-form ${\mathcal F}$ which restricts on each fibre of $Z\backslash D$  to the form $i\omega/s$ 
\item
the annihilator of ${\mathcal F}$ in $T_Z$ is the distribution generated by the vector field $Y$
\item
the residue is the  1-form $\varphi$ on $D$
\end{enumerate}
\end{prp}

\begin{prf} From the construction in Proposition \ref{algebroid} the transition functions $g_{UV}$ for $L_Z$ satisfy
$g_{UV}^{-1}dg_{UV}=(\varphi_V-\varphi_U)/s$ 
so $A_U=\varphi_U/i\zeta$ is a local connection form for a meromorphic connection on $L_Z$, with a simple pole on $\zeta=0$. There is a  similar form at $\zeta=\infty$ so the connection has a simple pole on the divisor $D$.  We see directly from the definition of $\varphi_U=i_{Y_U}\omega-f_Ud\zeta$ that the residue is the 1-form $\varphi$. 

The curvature form is ${\mathcal F}=dA_U=d(\varphi_U/s)$, but from (\ref{Fequation}) this can be expressed as $-i_Yd(\omega_U/s)$ and is thus annihilated by the vector field $Y$. If we establish Property 1, then ${\mathcal F}$ will be symplectic along the fibres, which are transversal to $Y$ and then we can deduce that $Y$ generates the annihilator.

Now $Y=Y_U+i\zeta Z_U$ and $\omega_U/s$ is closed on a fixed fibre $\zeta=c$, so $i_{Y_U}d(\omega_U/s)$ vanishes on the fibre since $Y_U$ by definition is tangential to a fibre.  We therefore only have to consider $-i\zeta i_{Z_U}d(\omega_U/s)$ for $\zeta=c$. Trivializing ${\mathcal O}(2)$ with $d/d\zeta$ we can write $s=i\zeta$, and using  $i_{Z{U}}\omega_U=0$ and ${\mathcal L}_{Z_U}\omega=0$ gives $i_{Z_U}d\omega_U=0$ modulo $d\zeta$. 
Hence on the fibre
$$-i\zeta i_{Z_U}d(\omega_U/s)=-i\zeta \omega_U i_{Z_U}d(1/i\zeta)=\omega_U/\zeta.$$
\end{prf}

\subsection{The link}
We shall now show that the holomorphic line bundle  with curvature $F=\omega_1+dd^c_1\mu$ on $Z=M\times S^2$ coincides with the holomorphic line bundle just constructed, trading \v Cech terminology for the Dolbeault viewpoint. To begin with, we  interpret in differential-geometric terms the Lie algebroid in Proposition \ref{alg} on each fibre of an arbitrary  hyperk\"ahler manifold. 

Without loss of generality consider $\zeta=0$ to be the complex structure with K\"ahler form $\omega_1$. On an open set $U$ we can write $2i\omega_1=\bar\partial\partial\phi_U=\bar\partial\theta_U$ for a $\partial$-closed $(1,0)$-form $\theta_U$. If $\omega_1$ is the curvature of a $U(1)$ connection, then $\theta_U$ is the connection form for a local holomorphic trivialization. On $U\cap V$, $\theta_U-\theta_V$ is therefore a $1$-cocycle of closed holomorphic $1$-forms and defines a holomorphic Lie algebroid. As connection forms $\theta_U-\theta_V=g_{UV}^{-1}\partial g_{UV}$ for transition functions $g_{UV}$.

\begin{prp} \label{algequiv} This is the Lie algebroid of Proposition \ref{alg} on the fibre at $\zeta=0$.
\end{prp}
\begin{prf} Given the hyperk\"ahler metric, the twistor space $Z$ is a $C^{\infty}$ product $M\times \CP^1$. We take the section $\omega$ of $\Lambda^2T_F^*(2)$ to be   
$$\omega=\left((\omega_2+i\omega_3)+2i\zeta\omega_1+\zeta^2(\omega_2-i\omega_3)\right)\frac{d}{d\zeta}$$
regarding the tangent vector $d/d\zeta$ as a local trivialization of ${\mathcal O}(2)$. 

Choose local holomorphic coordinates $z_1,\dots,z_n$ on the fibre $\zeta=0$. To fix terminology we shall call  $(1,0)$ vector fields  on a complex manifold those spanned  by $\partial/\partial z_i$.  Write $2i\omega_1=\bar\partial \theta_U$ as above and define the $(1,0)$ vector field $T_U$ on $U$ by $i_{T_U}(\omega_2+i\omega_3)=-\theta_U$, and put $X_k=[\partial/\partial\bar z_k, T_U]$. Since $\omega_2+i\omega_3$ is holomorphic,
\begin{equation}
i_{X_k}(\omega_2+i\omega_3)={\mathcal L}_{\partial/\partial\bar z_k}i_{T_U}(\omega_2+i\omega_3)= -{\mathcal L}_{\partial/\partial\bar z_k}\theta_U=    -2i(i_{\partial/\partial \bar z_k}\omega_1).
\label{Xone}
\end{equation}
The $(0,1)$ tangent vectors $X$ on $M$ for the complex structure  $\zeta$ are given by $i_X\omega=0$, and it therefore follows from (\ref{Xone}) that to first order $\partial/\partial \bar z_k+\zeta X_k$ is of type $(0,1)$. Furthermore, from the definition of $X_k$, 
$$\left[\frac{\partial}{\partial\bar z_k}+\zeta X_k,\frac{\partial}{\partial\zeta}+T_U\right]=0$$
at $\zeta=0$ so ${\partial}/{\partial\zeta}+T_U$ is a holomorphic vector field at $\zeta=0$ which projects to $d/d\zeta$ on $\CP^1$.
Thus, comparing with the construction in the previous section, we can take at $\zeta=0$  $Z_U={\partial}/{\partial\zeta}+T_U$.  

Now $\omega_2+i\omega_3$ is the symplectic form on the fibre, so 
 $$i_{X_{UV}}\omega=i_{Z_V-Z_U}\omega=i_{T_V-T_U}\omega=(\theta_U-\theta_V)$$
 showing we have have the same algebroid as in Proposition \ref{alg}.
\end{prf}

Now consider the case in question, where the hyperk\"ahler manifold has a circle action inducing a holomorphic vector field $Y$ on $Z$. 
As we have seen, the fibration  $\pi:Z\rightarrow \CP^1$ expresses the tangent bundle $T_Z$ as a holomorphic  extension
\begin{equation}
0\rightarrow T_F\rightarrow T_Z\rightarrow \pi^*T_{\CP^1}\rightarrow 0.
\label{extT}
\end{equation}
The product decomposition $Z=M\times \CP^1$ gives a $C^{\infty}$ decomposition of the tangent bundle as complex vector bundles $T_Z\cong T_M\oplus \pi^*T_{\CP^1}$ where $T_M$ at $(x,\zeta)$ is given the complex structure at $\zeta$. We denote by $\bar\partial_{\zeta}$ the $\bar\partial$-operator on $M$ with respect to the complex structure $\zeta$. Let $X^{1,0}$ denote the $(1,0)$ component with respect to this complex structure of the vector field $X$ (defining the circle action) on $M$. Then $i_X\omega=i_{X^{1,0}}\omega$ since $\omega$ is of type $(2,0)$ for all $\zeta$, and moreover since $\omega$ is symplectic on a fibre $i_{X^{1,0}}\omega$ uniquely defines $X^{1,0}$. 

For fixed $\zeta$, $\omega$ is holomorphic on $M$ so 
 $$i_{\bar\partial_{\zeta}X^{1,0}}\omega=\bar\partial_{\zeta}(i_{X^{1,0}}\omega)=(di_X\omega)^{1,1}.$$
But, as a form on the fibre $M$,
$$di_X\omega={\mathcal L}_X\omega={\mathcal L}_X\left((\omega_2+i\omega_3)+2i\zeta\omega_1+\zeta^2(\omega_2-i\omega_3)\right)=i(\omega_2+i\omega_3)-i\zeta^2(\omega_2-i\omega_3)$$
using Equation (\ref{Xomega123}) for the action on $\omega_1,\omega_2,\omega_3$. 
Since $\omega$ is of type $(2,0)$, 
$(\omega_2+i\omega_3)^{1,1}=(-2i\zeta\omega_1-\zeta^2(\omega_2-i\omega_3))^{1,1}$. It follows that  for $\zeta\ne \infty$ 
$i_{\bar\partial_{\zeta}X^{1,0}}\omega=(di_X\omega)^{1,1}$
 is divisible by $\zeta$ and hence so is  $\bar\partial_{\zeta}X^{1,0}$. A similar argument near $\zeta=\infty$ (or using the real structure) shows that $\zeta^{-1}\bar\partial_{\zeta}X^{1,0}$ is a well defined $(0,1)$-form on $Z$ with values in $T_F(-2)$.

These facts play a useful role because of the following description of the tangent bundle of $Z$, the holomorphic extension $0\rightarrow T_F\rightarrow T_Z\rightarrow \pi^*T_{\CP^1}\rightarrow 0$, as a $C^{\infty}$ direct sum. 
\begin{lemma}\label{holstr} The holomorphic structure of $T_Z\cong T_F\oplus \pi^*T_{\CP^1}$ is  defined by the operator $\bar\partial_Z$ where 
$$\bar\partial_Z(W,u)=(\bar\partial_{\zeta}W+iu{\zeta}^{-1}\bar\partial_{\zeta}X^{1,0}, \bar\partial_{\zeta} u).$$
\end{lemma}
\begin{lemprf}
The  vector field  $Y$ gives a holomorphic splitting of  the extension  (\ref{extT}) outside the divisor $D=Z_0+Z_{\infty}$ since $Y$ projects to the non-zero tangent vector $i\zeta d/d\zeta$ on $\CP^1$. Writing $s$ as the corresponding section of ${\mathcal O}(2)$, this means that the extension is defined by the class  $\delta(Y_D)\in H^1(Z,T_F(-2))$ in the long exact cohomology sequence of 
\begin{equation}
0\rightarrow T_F(-2)\stackrel{s}\rightarrow T_F\rightarrow T_F\vert_D\rightarrow 0
\label{ext2}
\end{equation}
where $Y_D$ is the vector field $Y$ on $D$, where it is tangential to those two fibres. 

On the other hand, since $X^{1,0}$ is holomorphic in $\zeta$,  $-i\zeta^{-1}\bar\partial_{\zeta}X^{1,0}=-i\zeta^{-1}\bar\partial_Z X^{1,0}=\gamma$  is a Dolbeault representative for a class in $H^1(Z,T_F(-2))$. Then $s\gamma=\bar\partial_ZX^{1,0}$ is cohomologically trivial and in the Dolbeault version of the exact cohomology sequence is   defined by $\delta(X^{1,0}\vert_D)$. But on $D$, $X$ is the holomorphic vector field $Y_D$,  
so the form $\gamma$ defines the extension and the lemma follows.  
\end{lemprf}
On the cotangent bundle this implies that 
the $\bar\partial$-operator is 
\begin{equation}
\bar\partial_Z( \alpha, u)= (\bar\partial_{\zeta} \alpha, \bar\partial_{\zeta}u+i{\zeta}^{-1}\alpha(\bar\partial_{\zeta}X^{1,0})).
\label{cot}
\end{equation}

To prove the equivalence between the differential geometric approach and the twistorial approach, we need to show that the holomorphic Lie algebroid on $Z$ defined by the closed $(1,1)$-form $F=\omega_1+dd_1^c\mu$ is obtained by a coboundary map from a holomorphic section of $T_Z^*(-2)$ on $D$.

Observe first that for any $1$-form $\alpha$ on $M$, $(2i\zeta I+(1+\zeta^2)J+i(1-\zeta^2)K)\alpha$ is of type $(1,0)$ with respect to the complex structure $\zeta$. Take $\alpha=d\mu$ and define the $(1,0)$ form $\phi$ on $M$ by
\begin{equation}
\phi=2i\zeta d_1^c\mu+(1+\zeta^2)d_2^c\mu+i(1-\zeta^2)d_3^c\mu.
\label{fidef}
\end{equation}
and use the $C^{\infty}$ splitting of $T_Z^*$ to interpret this in $\Omega^{1,0}(Z,{\mathcal O}(2))$ (the twist by ${\mathcal O}(2)$ comes from the quadratic dependence on $\zeta$). 

Now define a section $\psi$ of $T_Z^*(2)\cong T^*_F(2)\oplus \pi^*T^*_{\CP^1}(2)$ by $\psi=(\phi,-2i\mu)$. 

\begin{lemma} $\bar\partial_Z\psi=2i\zeta F$
\end{lemma}
\begin{lemprf}
 First we check the result  restricted to fibres: i.e. that $\bar\partial_{\zeta}\phi=2i\zeta F$. Using $dd^c_3\mu=-\omega_3$ and  $dd^c_2\mu=-\omega_2$, we see that  $\bar\partial_{\zeta}\phi$ is the $(1,1)$-component of 
$$d\phi=2i\zeta dd_1^c\mu-(1+\zeta^2)\omega_2-i(1-\zeta^2)\omega_3.$$
But $(\omega_2+i\omega_3)+2i\zeta\omega_1+\zeta^2(\omega_2-i\omega_3)$ is of type $(2,0)$ 
 so, as before 
 $$((1+\zeta^2)\omega_2+i(1-\zeta^2)\omega_3)^{1,1}=-2i\zeta\omega_1^{1,1}$$
hence 
 \begin{equation}
 \bar\partial_{\zeta}\phi  = 2i\zeta( dd_1^c\mu+\omega_1)^{1,1}=2i\zeta F
 \label{partfi}
 \end{equation} 
 since $F$ is of type $(1,1)$ with respect to {\it all} complex structures. 
 
 To finish we have to show, using   (\ref{cot}), that $2\bar\partial_{\zeta}\mu=\zeta^{-1}\phi(\bar\partial_{\zeta}X^{1,0})$.
 Now, contracting the 1-form $\phi$ with  $\bar\partial_{\zeta}X^{1,0}\in \Omega^{0,1}(T)$ we have, using (\ref{partfi}), 
 
\begin{equation}
\phi(\bar\partial_{\zeta}X^{1,0})=\bar\partial_{\zeta}(i_X\phi)+i_{X^{1,0}}\bar\partial_{\zeta}\phi=\bar\partial_{\zeta}(i_X\phi)+2i\zeta(i_XF)^{0,1}
\label{fiagain}
\end{equation}  
 and by definition 
$$i_X\phi=-2i\zeta g(X,X)-(1+\zeta^2)g(X,JIX)+i(1-\zeta^2)g(X,KIX)=-2i\zeta g(X,X).$$
But $i_XF=d(\mu+g(X,X))$ so $(i_XF)^{0,1}=\bar\partial_{\zeta}(\mu+g(X,X))$ and so from (\ref{fiagain}) 
$$\phi(\bar\partial_{\zeta}X^{1,0})=-2i\zeta \bar\partial_{\zeta}g(X,X) +2i\zeta \bar\partial_{\zeta}(\mu+g(X,X))=2i\zeta \bar\partial_{\zeta}\mu$$
 \end{lemprf}
 
 From the Lemma we see that the algebroid given by $F$ is defined by $\delta(\psi_D)/2$ where $\psi_D$ is the holomorphic section of $T_Z^*(2)$ on $D$ given by the restriction of $\psi$. Now at $\zeta=0$, we see from (\ref{fidef}) that $\phi$ restricts to 
 $$d_2^c\mu+id_3^c\mu=-i_X\omega_2+i(i_X\omega_3)=-i(i_X(\omega_2+i\omega_3))$$
 from the properties of $\mu$ in Section \ref{diff}. This is the section $i_Y\omega$  of $T_F^*(2)$ on $Z_0$ in the twistor construction, and the choice of moment map $\mu$ is equivalent to the choice of lift of $Y$ to the algebroid.
  
 Finally the integrality condition on the cohomology class of $\omega_1$ is clearly the same as for $F=\omega_1+dd^c_1\mu$ and so if satisfied, the algebroid on $Z$ defines a holomorphic principal $\C^*$-bundle.  
 
 To summarize:
 \begin{prp} The line bundle $L_Z$ on $Z$ which gives the hyperholomorphic connection with curvature $\omega_1+dd^c_1\mu$ on $M$ is obtained by applying the construction of Section \ref{twist} to the section $\omega=\left((\omega_2+i\omega_3)/2+i\zeta\omega_1+\zeta^2(\omega_2-i\omega_3)/2\right){d}/{d\zeta}$ of $\Lambda^2T^*_F(2)$.
  \end{prp} 
 \section{Examples}
 \subsection{The flat case}\label{flat}
 We take $M$  to be $\C^{2k}$ -- a product of $k$ copies of the motivating example in Section \ref{diff}. So the hyperholomorphic line bundle  has curvature 
 $$F=\frac{i}{2}\sum_i( dz_i\wedge d\bar z_i-dw_i\wedge d\bar w_i).$$
 The twistor space $Z$ is then the total space of the vector bundle $\C^{2k}(1)$ over $\CP^1$. We write this more invariantly as $Z=W(1)\oplus W^*(1)$ where $W$ is a $k$-dimensional vector space. Then the natural pairing $\langle v,\xi\rangle$ defines the section $\omega$ of $\Lambda^2T^*_F(2)$, which is just a constant skew form on the fibres, and the $S^1$-action extends to a $\C^*$-action which is the composition of the natural action on $\CP^1$ with the action $(v,\xi)\mapsto (v,\lambda \xi)$. A hermitian structure on $W$ defines an isomorphism $W\cong \bar W^*$ which, together with the antipodal map on $\CP^1$ generates the real structure on $Z$. We shall now describe the holomorphic line bundle on $Z$ constructed in  Section \ref{twist}, and to do this in \v Cech language  we describe $Z$ in terms of two open sets.
 
 Define $U,V\subset Z$ to be the subsets $\zeta\ne \infty$ and $\zeta\ne 0$ respectively. Then we have coordinates $v_i,\xi_i,\zeta$ on $U$ and 
  $\tilde v_i,\tilde \xi_i,\tilde\zeta$ on $V$ where on $U\cap V$,
  $$\tilde\zeta=1/\zeta,\qquad \tilde v_i=v_i/\zeta,\qquad \tilde\xi_i=\xi_i/\zeta.$$
  The real structure is $(v_i,\xi_i,\zeta)\mapsto (\bar\xi_i/\bar\zeta,-\bar v_i/\bar\zeta,-1/\bar\zeta)$ 
  and, comparing with Example \ref{ex1},  $\omega_U=\sum_idv_i\wedge d\xi_i/2$ and similarly over $V$.
  
 The $\C^*$-action is given by $(v_i,\xi_i,\zeta)\mapsto (v_i,\lambda\xi_i,\lambda\zeta)$ for $\zeta\ne \infty$ and correspondingly  $(\tilde v_i,\tilde\xi_i,\tilde\zeta)\mapsto (\lambda^{-1}\tilde v_i,\tilde\xi_i,\lambda^{-1}\tilde\zeta)$.
 Over $U$  and $V$ we have the vector field $Y$ on $Z$ expressed as 
 $$Y=\sum _i i\xi_i\frac{\partial}{\partial \xi_i}+i\zeta\frac{\partial}{\partial \zeta}=-\sum _i i\tilde v_i\frac{\partial}{\partial \tilde v_i}-i\tilde\zeta\frac{\partial}{\partial \tilde\zeta}.$$
We take  $Z_U=\partial/\partial \zeta$ and $Z_V=\partial/\partial \tilde\zeta$ to be the lifted vector fields. Since  $\omega_U$ in these coordinates is independent of $\zeta$ this is the required condition. Over $U\cap V$, the cocycle $X_{UV}$ with values in $T_F(-2)$ is 
$$X_{UV}=\left( v_i\frac{\partial}{\partial v_i}+ \xi_i\frac{\partial}{\partial \xi_i}\right)\frac{d\zeta}{\zeta}.$$
The open sets $U,V$ are $\C^*$-invariant, as is $X_{UV}$, so we may take the $f_U=0$ and then the cocycle which defines the Lie algebroid on $Z$ is  
 \begin{equation} \frac{1}{2i\zeta}(i_{Y_U}\omega_U-i_{Y_V}\omega_V)=-\frac{1}{2\zeta}\sum_i(\xi_idv_i+\zeta^2\tilde v_id\tilde\xi_i).
 \label{start}
 \end{equation}
Note that the factor $\zeta^2$ appears to relate the two local trivializations of ${\mathcal O}(2)$ over $U$ and $V$.

 This expands in the $U$-coordinates to 
$$- \frac{1}{2\zeta}\sum_i(\xi_idv_i+\zeta v_id(\xi_i/\zeta))=-d\sum_i v_i\xi_i/2\zeta.$$
Thus the line bundle is defined by the transition function $g_{UV}$ where 
$$g_{UV}=\exp(-\sum_i v_i\xi_i/2\zeta).$$

The meromorphic connection is given by one-forms $A_U,A_V$ such that 
$$A_V-A_U=g_{UV}^{-1}dg_{UV}=-d\sum_iv_i\xi_i/2\zeta.$$
But this can be read off from (\ref{start})
$$A_U=\sum_i \xi_idv_i/2\zeta, \quad A_V=-\sum_i\tilde v_id\tilde \xi_i/2\tilde\zeta.$$

That this describes the hyperholomorphic line bundle is a consequence of the general proof in the previous section, but it is instructive to see it directly. The $C^{\infty}$ description of the twistor space as a product is given by setting
$$v_i=z_i+\zeta \bar w_i,\quad \xi_i=w_i-\zeta \bar z_i$$
for these equations define the real holomorphic sections of $\C^{2k}(1)$ parametrized by coordinates $(z_i,w_i)\in \C^{2k}$. A hermitian metric on a line bundle with holomorphic transition functions $g_{UV}$ is give by smooth real-valued functions $h_U$ on $U$ such that on $U\cap V$, $ h_U=\vert g_{UV}\vert^2 h_V$. In our case
$$\log \vert g_{UV}\vert =-\frac{1}{2}\Real \sum_i v_i\xi_i/\zeta=-\frac{1}{2}\Real \sum_i z_iw_i/\zeta+w_i\bar w_i-z_i\bar z_i-\zeta\bar z_i\bar w_i$$
so we can take
\begin{equation}
\log h_U=\frac{1}{2}\sum_i z_i\bar z_i-w_i\bar w_i+\zeta\bar z_i\bar w_i+\bar\zeta z_i w_i
\label{herm}
\end{equation}
and $\log h_V=-\log h_U(-1/\bar\zeta)$.

Now the $(1,0)$-forms on $Z$ for $\zeta\ne \infty$ are spanned by $dz_i+\zeta d\bar w_i,dw_i-\zeta d\bar z_i, d\zeta$ and a short calculation gives 
$$\bar\partial_Z \log h_U= \frac{1}{2}\sum z_iw_id\bar\zeta+z_id\bar z_i-w_id\bar w_i+\bar\zeta d(z_iw_i)$$
and hence 
$\bar\partial_Z \partial_Z \log h_U=(\sum _i-dz_id\bar z_i+dw_id\bar w_i)/2$, giving the curvature of the hyperholomorphic connection on the line bundle pulled back to the twistor space. 
\subsection{Taub-NUT space}
There is another well-known hyperk\"ahler metric on $\R^4$ whose twistor space can be described easily, namely the Taub-NUT metric.  It is  the submanifold of the vector bundle $L(1)\oplus L^*(1)$ over ${\mathcal O}(2)$ defined by  $\{(v,\xi): v\xi=\eta\}$ where $\eta$ is the tautological section of ${\mathcal O}(2)$ pulled back to its total space, and $L$ is the line bundle with transition function $\exp(\eta/\zeta)$ (the line bundle occurring in the construction of monopoles \cite{Hit1}). It has local coordinates $v,\xi,\zeta$ and $\tilde v,\tilde \xi,\tilde\zeta$ as above where
$$\tilde\zeta=1/\zeta;\quad \tilde v=\exp(v\xi/\zeta)v/\zeta;\quad \tilde \xi=\exp(-v\xi/\zeta)\xi/\zeta.$$
The $\C^*$-action is as in the flat case and a similar calculation to the one above gives  the cocycle
$$- \frac{1}{\zeta}\left(\xi dv+\zeta v[d(\xi/\zeta)-d(v\xi/\zeta)\xi/\zeta]\right)=d(-v\xi/\zeta+(v\xi/\zeta)^2)).$$
The hyperholomorphic line bundle is therefore given by the transition function 
$$g_{UV}=\exp(-\eta/\zeta+(\eta/\zeta)^2)).$$

  \subsection{The semi-flat case}
  We next describe from a differential geometric point of view the line bundle as it appears in the physicist's c-map construction (see \cite{cort} for further discussion of this). This starts with a {\it  special K\"ahler manifold}. This is a manifold ${\mathcal M}$ with a flat symplectic connection $\nabla$ with symplectic form $\omega$, and a Hamiltonian vector field $X$ such that $\nabla X=I$ is a complex structure making $\omega$ a K\"ahler form.  When $X$ preserves the metric, or equivalently the complex structure, and generates a circle action, it is called a conical special K\"ahler manifold and  its symplectic quotient is called a projective special K\"ahler manifold. We shall assume that the metric is  positive definite, though in many of the moduli space occurrences of this structure it is indefinite.

  Let $\phi$ be the Hamiltonian function generating $X$, then following the description in \cite{Hit2} we let $x_1,\dots, x_{2k}$ be real flat coordinates for the connection $\nabla$. In these coordinates  the symplectic  form is given by a constant matrix $\omega_{ij}$. On the product with $\R^{2k}$, with coordinates $y_1,\dots,y_{2k}$ we define on ${\mathcal M}\times \R^{2k}$ three closed $2$-forms by 
  $$\omega_1+i\omega_2=\sum_{jk}\omega_{jk}d(x_j+iy_j)\wedge d(x_k+iy_k)\quad \omega_3=\sum_{jk}\frac{\partial^2\phi}{\partial x_j\partial x_k}dx_j\wedge dy_k.$$
  Then these define a hyperk\"ahler metric. Take the trivial action of the circle on the $\R^{2k}$ factor and this gives an action which  leaves fixed $\omega_1$, and which has moment map $\mu=\phi$, a function of ${\mathcal M}$ alone.  
  
  The metric on ${\mathcal M}$ is of Hessian form $g=\nabla^2\phi$, so since $Id\phi=\grad\phi$ we have 
  $$Id\phi=\sum_i X_ig_{ij}dx_j=\sum_{ijk}\omega^{ik}\frac{\partial\phi}{\partial x_k}\frac{\partial^2\phi}{\partial x_i\partial x_j}dx_j$$
  and so
  $$dId\phi=\sum_{ijk}\omega^{ik}\frac{\partial^2\phi}{\partial x_{\ell}\partial x_k}\frac{\partial^2\phi}{\partial x_i\partial x_j}dx_{\ell}\wedge dx_j=\sum_{ij}\omega_{ij}dx_i\wedge dx_j$$
 since $\omega$ is a K\"ahler form for the  Hessian  metric.  
 Hence
 $$\omega_1+dd^c_1\mu=\sum_{jk}\omega_{jk}(dx_j\wedge dx_k-dy_j\wedge dy_k)+dd^c_1\mu=-\sum_{jk}\omega_{jk}dy_j\wedge dy_k.$$
  The curvature form is thus constant on each copy of $\R^{2k}$.
  
  In many situations, there is a lattice in the fibres and the hyperk\"ahler metric is defined on the quotient which is a torus fibration over ${\mathcal M}$ -- in fact a holomorphic integrable system for the symplectic form $\omega_2+i\omega_3$. The hyperholomorphic line bundle then defines a complex line bundle over each torus. The corresponding principal bundle is a Heisenberg extension of the torus group.
  
  The metric above has a $2k$-dimensional  abelian group of translational triholomorphic symmetries (preserving $I,J$ and $K$) generated by $\partial/\partial y_i$. It is a special case of the Legendre transform construction which requires only a $k$-dimensional group, and we shall give a twistor description in the next section.
   \subsection{The Legendre transform}\label{Leg}
   The Legendre transform \cite{HKLR} is a means of constructing a hyperk\"ahler metric on a manifold $M^{4k}$ with a $k$-dimensional abelian group $G$ of triholomorphic symmetries. It requires a hyperk\"ahler moment map for the group which  means that the restriction of the K\"ahler forms $\omega_1,\omega_2,\omega_3$ to the orbits must be zero. In particular to apply it to the previous semi-flat case, we need to restrict the translations to a Lagrangian subspace of $\R^{2k}$ with respect to the symplectic form. 
   
   If  $G=W\cong \R^k$ and acts freely then the hyperk\"ahler moment map expresses $M$ as a principal $W$-bundle over $W^*\otimes \R^3$. The twistor space  is then a principal $\C^k$-bundle over the vector bundle $W^*(2)\cong \C^k(2)\rightarrow \CP^1$. The construction starts with a class  $a\in H^1(W^*(2),{\mathcal O(2)})$, then differentiating along the fibres we get a class $d_Fa\in W\otimes H^1(W^*(2),{\mathcal O})$. This class defines the principal $\C^k$-bundle. 
   
   We take again the  open sets $U,V$ on $W^*(2)$ defined by $\zeta\ne \infty,\zeta\ne 0$ and coordinates $\eta_i,\zeta$ on $U$ and $\tilde\eta_i,\tilde\zeta$ on $V$ and on $U\cap V$ we have 
   $$\tilde\eta_i=\eta_i/\zeta^2,\qquad \tilde\zeta=1/\zeta$$
   
   \begin{rmk} Most of the interesting examples do not quite fall into this picture because the group action has fixed points: the cohomological formulation is replaced by contour integrals and these may be of multi-valued functions, but for our purposes we shall stay with this global picture.
   \end{rmk}
    
    The class $a\in H^1(W^*(2),{\mathcal O(2)})$ is defined by a cocycle on $U\cap V$ of the form
    $$H(\eta_1,\dots,\eta_k,\zeta)\frac{d}{d\zeta}$$ and then the principal $\C^k$-bundle has extra coordinates $\chi_1,\dots,\chi_k$ over $U$  and $\tilde\chi_1,\dots,\tilde \chi_k$ over $V$ which on $U\cap V$ satisfy
    $$\tilde\chi_i=\chi_i+\frac{\partial H}{\partial \eta_i}$$
    This preserves the symplectic form $\omega=\sum_id\chi_i\wedge d\eta_i=\zeta^2\sum_id\tilde\chi_i\wedge d\tilde\eta_i$ along the fibres.

    \begin{ex} Take $M=\C^2$  and let the action of $W\cong\R$ be $(z,w)\mapsto (z+t,w)$. Then on the twistor 
    space $Z=\C^2(1)$ the $\C$-action in coordinates is given by $(v,\xi)\mapsto (v+\zeta,\xi+1)$. Then $\eta=\zeta\xi-v$ is invariant and defines the projection from $\C^2(1)$ to ${\mathcal O}(2)$ over $U$ and $\tilde\eta=\tilde\xi-\tilde\zeta\tilde v$ over $V$. Over $U$ we have a section $(v,\xi)=(-\eta,0)$ and over $V$, $(\tilde v,\tilde \xi)=(0,\tilde\eta)$. Then on $U\cap V$, $(0,\tilde\eta)=(-\eta,0)+\chi(\zeta,1)$ 
    giving $\chi=\eta/\zeta$ and hence $H(\eta,\zeta)=\eta^2/2\zeta$.
    \end{ex}
    
   To determine the holomorphic line bundle over the twistor space, we take the natural lift of the $\C^*$-action on $\CP^1$ to  $W^*(2)$ which in coordinates is $\zeta\mapsto\lambda\zeta,\eta_i\mapsto \lambda \eta_i$. Since $U$ and $V$ are invariant we can average a cocycle over the circle group to make it invariant. Then $H(\eta_1,\dots,\eta_k,\zeta)$ is homogeneous of degree $1$ and as before we can take the $f_U=0$. Now we have 
    $$\frac{1}{i\zeta}(i_{Y_U}\omega_U-i_{Y_V}\omega_V)=-\frac{1}{\zeta}\sum_i(\eta_id\chi_i-\zeta^2\tilde \eta_id\tilde\chi_i)=\frac{1}{\zeta}\sum_i\eta_id \left(\frac{\partial H}{\partial \eta_i}\right).$$
    \begin{prp} The hyperholomorphic line bundle is defined by the transition function 
    $$g_{UV}=\exp\left(-\frac{1}{2}\frac{\partial H}{\partial \zeta}\right).$$
   \end{prp}
   \begin{prf}
   $$\frac{\eta_i}{\zeta}d\left(\frac{\partial H}{\partial \eta_i}\right)=d\left(\frac{\eta_i}{\zeta}\frac{\partial H}{\partial \eta_i}\right)-\frac{\partial H}{\partial \eta_i}\frac{d\eta_i}{\zeta}+\eta_i\frac{\partial H}{\partial \eta_i}\frac{d\zeta}{\zeta^2}$$
   But from the homogeneity of $H$,
   $$\sum_i {\eta_i}\frac{\partial H}{\partial \eta_i}+\zeta\frac{\partial H}{\partial \zeta}=H.$$
   Substituting, we get
     \begin{eqnarray*}\sum_i\frac{\eta_i}{\zeta}d\left(\frac{\partial H}{\partial \eta_i}\right)&=&d\left(\frac{H}{\zeta}-\frac{\partial H}{\partial \zeta}\right)-\frac{1}{\zeta}dH+\frac{1}{\zeta}\frac{\partial H}{\partial \zeta}d\zeta+H\frac{d\zeta}{\zeta^2}-\frac{1}{\zeta}\frac{\partial H}{\partial \zeta}d\zeta\\
     &=&d\left(-\frac{\partial H}{\partial \zeta}\right)
       \end{eqnarray*}
   \end{prf}
   \begin{ex} In the flat example above ${\partial H}/\partial \zeta=-\eta^2/2\zeta^2=-(\zeta\xi-v)^2/2\zeta^2$ which expands to  $-\xi^2/2+v\xi/\zeta-v^2/2\zeta^2$. But $\exp \xi^2/2$ is a gauge transformation over $U$ and $\exp(v^2/2\zeta^2)=\exp(\tilde v^2/2)$ is well defined over $V$, so after changing local trivializations we get transition function $\exp(-v\xi/2\zeta)$, the same holomorphic bundle as  in Section \ref{flat}. 
   \end{ex}
   
   In the semi-flat case with a circle action we have a bigger group of symmetries.  On the space $W^*(2)$ these are realized by transformations $\eta_i\mapsto \eta_i+t_i\zeta$. We cannot impose the symmetry on the cocycle itself since $\zeta\partial H/\partial\eta_i=0$ is too strong, but we can ask that  $\zeta H$ (and hence also $\partial H/\partial\eta_i$) must be  a coboundary. This means that it comes from a class in $H^0(D,{\mathcal O}(4))$ in the exact sequence $$\rightarrow H^0(W^*(2),{\mathcal O}(4)) \rightarrow H^0(D,{\mathcal O}(4))\rightarrow H^1(W^*(2),{\mathcal O}(2))\rightarrow$$ where $D$ is the sum of the divisors $\zeta=0$ and $\zeta=\infty$. The cocycle is then of the form 
   $$H=\frac{1}{\zeta}\left(f(\eta_1,\dots,\eta_k)+\bar f(\tilde\eta_1,\dots,\tilde\eta_k))\right)$$
    for a holomorphic function $f$, homogeneous of degree $2$.  (In fact, (see \cite{Hit2}) this is essentially the holomorphic prepotential of special K\"ahler geometry.) Then the transition function is 
   $$\exp\left(-\frac{1}{2}\frac{\partial H}{\partial \zeta}\right)=\exp\left(\frac{1}{2\zeta^2}(f(\eta_1,\dots,\eta_k)+\bar f(\tilde\eta_1,\dots,\tilde\eta_k))\right).$$
  
    \subsection{Cotangent bundles}\label{cotangent}
    The cotangent bundle of a K\"ahler manifold $N$ is a holomorphic symplectic manifold with a circle action along the fibres, which acts by multiplication on the symplectic form, as in the cases here. Feix \cite{Feix1}, and independently Kaledin and Verbitsky, showed that if the metric on $N$ is real analytic then there is a canonical extension to an $S^1$-invariant  hyperk\"ahler metric  on a neighbourhood of the zero section of $T^*N$. Any hyperk\"ahler metric is real analytic so this is a necessary condition anyway. 
    Since this is local around the zero section the same argument also applies to hyperk\"ahler metrics with a circle action with fixed point set which has the local structure of the cotangent bundle \cite{Feix2}. In particular the hyperk\"ahler metric is uniquely determined by the  metric on the fixed point set. 
    
    In a further paper \cite{Feix3}, Feix showed that a real analytic connection on a bundle over a K\"ahler manifold $N$, whose curvature is of type $(1,1)$, extends uniquely to an $S^1$-invariant  hyperholomorphic bundle on this hyperk\"ahler extension on $T^*N$. 
    
    Now the hyperholomorphic bundle considered here has curvature $F=\omega_1+dd^c_1\mu$. The function $\mu$ is the moment map for the circle action and so is constant on the fixed point set $N$. The action also preserves the complex structure $I$ so the fixed point set is complex and hence  $dd^c_1\mu$ restricted as a form vanishes there. It follows that the hyperholomorphic bundle on $T^*N$ is the unique extension of the line bundle with curvature the K\"ahler form on $N$.  Since Feix's construction uses twistor theory, we can implement our description in this rather general case too. We first outline Feix's argument. 
    
    Locally, real analyticity allows a complexification $N^c$ of the K\"ahler manifold $N$. This means that $z_i,\bar z_i$ become independent coordinates 
    $z_i,\tilde z_i$ and there are two transverse  foliations whose leaves are defined by $z_i=c_i$ and $\tilde z_i=\tilde c_i$. Let $B_+$ be the quotient space of $N^c$ by the first and $B_-$ of the second. Then the projections $p_{\pm}:N^c\rightarrow B_{\pm}$ identify the real submanifold $N\subset N^c$ with $B_{\pm}$ and give it the complex structures $\pm I$. 
    
    The K\"ahler form is now a holomorphic symplectic form on $N^c$ and the foliations are Lagrangian foliations. This means the fibres have natural flat affine structures: functions on the base $B_{\pm}$ define Hamiltonian vector fields along the fibres and these are covariant constant with respect to the flat affine connection. Let $V_{\pm}$ denote the rank $k+1$ vector bundle over $B_{\pm}$ consisting of the affine linear functions along the fibres. In canonical coordinates, if $\omega=\sum_i dq_i\wedge dz_i$ then this space is spanned  by $\{1,q_1,\dots,q_k\}$. 
    
    There is a natural map $\phi_+:N^c\times \C^*\rightarrow V_+^*$ defined by $\phi_+(x,\lambda)=\lambda\delta_x$ where $\delta_x(f)=f(x)$ for an affine linear function $f$. This is locally biholomorphic. We define $\phi_-$ similarly and then construct a twistor space $Z$ by attaching $V_+^*$ to $V_-^*$ by identifying $\phi_+(x,\lambda)$ with $\phi_-(x,\lambda^{-1})$ over a neighbourhood of $S^1$ in $\C^*$. Note that the constant functions are contained in $V_{\pm}$ which provides  homomorphisms $V^*_{\pm}\rightarrow \C$. The identification then gives the projection   from $Z$ to $\CP^1$. 
    
    Now take a real  analytic  K\"ahler potential $f(z,\bar z)$, so that the K\"ahler form satisfies $2i\omega=\bar\partial\partial f $. In the complexification this can be written as $2i\omega=\sum_j  d(\partial f/\partial z_j)\wedge dz_j$ and then the affine linear functions on a leaf $z_i=c_i$ are linear combinations of $1$ and $q_i=(\partial f/\partial z_i)(c,\tilde z)$ for $i=1,\dots,k$. With respect to this basis we can take local coordinates $a_0,a_1,\dots, a_k, z_1,\dots z_k$ on $V_+^*$ and similarly $\tilde a_i,\tilde z_j$ on $V_-^*$. Then the identification to define $Z$ can be written  as 
    $$a_0=\zeta, \quad a_i=\zeta\frac{\partial f}{\partial z_i},\qquad \tilde a_0=\zeta^{-1},\quad \tilde a_i=\zeta^{-1}\frac{\partial f}{\partial \tilde z_i}.$$
    In other words, we have coordinates $a_0,a_1,\dots, a_k, z_1,\dots z_k$ on an open set $U$ and $\tilde a_0,\dots, \tilde a_k, $ $\tilde z_1,\dots \tilde z_k$ on $V$ where $a_i=\zeta{\partial f}/{\partial z_i}$ defines $\tilde z_i$ as a function of $a_0,a_1,\dots, a_k, z_1,\dots z_k$.
    
    The further properties to recognize this as the twistor space of a hyperk\"ahler manifold are attended to in \cite{Feix1}, but for us we only need to note  
  the  natural $\C^*$-action by multiplication on the space of affine linear functions.

 \begin{prp} Let $f$ be a K\"ahler potential on the real analytic K\"ahler manifold $(N,\omega)$, then a transition function for the line bundle $L_Z$ on the twistor space of $T^*N$ which defines the hyperholomorphic extension of the line bundle on $N$ with curvature $\omega$ is $$g_{UV}=\exp (-f(z,\tilde z)/2)$$ where $\tilde z_i$ is expressed as a function of $a_0,a_1,\dots, a_k, z_1,\dots z_k$ by $a_i=\zeta{\partial f}/{\partial z_i}$.
 \end{prp}
 
 \begin{prf} 
 From the construction above, and using the terminology of Section \ref{twist}  we have $i_{Y_U}\omega_U=-i\sum \zeta{\partial f}/{\partial z_i}dz_i$ and $i_{Y_V}\omega_V=i\sum \zeta^{-1}{\partial f}/{\partial \tilde z_j}d\tilde z_j$. Hence the cocycle is 
 $$\frac{1}{2i\zeta}(i_{Y_U}\omega_U-i_{Y_V}\omega_V)=-\frac{1}{2}\left(\frac{1}{\zeta}\sum_i  \zeta\frac{\partial f}{\partial z_i}dz_i-\zeta^2 \zeta^{-1}\frac{\partial f}{\partial \tilde z_i}d\tilde z_i\right)=-\frac{1}{2}df$$
  \end{prf} 
  
  We may compare this formula with the hyperholomorphic extension of the connection on a line bundle $L$ with curvature $\omega$. The argument of    \cite{Feix3} is as follows. Since $\omega$ is of type $(1,1)$ the complexification vanishes on restricting the form to the leaves of each foliation, and in simply-connected neighbourhoods there exist line bundles $L_{\pm}$ on $B_{\pm}$ consisting of sections of $L$ on $N^c$ covariant constant  along the leaves.  If $q_{\pm}:V^*_{\pm}\rightarrow B_{\pm}$ is the projection then a line bundle on the twistor space is obtained by identifying $q_+^*L_+$ and $q_-^*L_-$ as follows. For $x\in N^c$, $\phi_+(x,\lambda)$ is identified with $\phi_-(x,\lambda^{-1})$ and a vector in $(q_+^*L_+)_{\phi_+(x,\lambda)}$, which is a covariant constant section $\sigma_+$, is identified with a covariant constant section $\sigma_-$ if $\sigma_+(x)=\sigma_-(x)$. 
  
  Now 
  $$2i\omega=-\sum_{i,j}  \frac{\partial^2 \!f}{\partial z_i\partial \tilde z_j}dz_i\wedge d\tilde z_j=-\sum_{j}  d(\frac{\partial \!f}{\partial \tilde z_j} d\tilde z_j)$$
    so that $-\sum_{j} ( {\partial \!f}/{\partial \tilde z_j}) d\tilde z_j/2$ is a connection form for the line bundle on the leaves $z_i=c_i$, and $\exp( f(c,\tilde z)/2)$ is a covariant constant section.  The identification is therefore achieved by $\exp(- f(z,\tilde z)/2)$ as in the proposition.

    \begin{ex} Flat space $\C^{2k}$ has $f=\sum_iz_i\bar z_i$ complexifiying to $f=\sum_iz_i\tilde z_i$. So 
    $$a_0=\zeta, \quad a_i=a_i=\zeta\frac{\partial f}{\partial z_i}=\zeta \tilde z_i\qquad \tilde a_0=\zeta^{-1},\quad \tilde a_i=\zeta^{-1}\frac{\partial f}{\partial \tilde z_i}=\zeta^{-1} z_i$$ 
    and the twistor space has coordinates $a_0,a_1,\dots,a_k,z_1,\dots, z_k$ on $U$ and
    $$ \tilde a_0=\zeta^{-1}\quad \tilde a_i=\zeta^{-1} z_i\quad  \tilde z_i=\zeta^{-1}a_i$$
    on $V$ giving $\C^{2k}(1)\rightarrow \CP^1$. The transition function for the hyperholomorphic line bundle is $\exp (-\sum z_ia_i/2\zeta)$ as in the direct calculation in Section \ref{flat}. 
    \end{ex}
     \subsection{Monopole moduli spaces}
     Among the gauge-theoretic  hyperk\"ahler moduli spaces, the one we know most about is the moduli space of $SU(2)$ magnetic monopoles of charge $k$ on $\R^3$ \cite{AH}. These are solutions to the Bogomolny equations $F=\ast\nabla\phi$ where $F$ is the curvature of a connection  $A$ and $\phi$, the Higgs field, is a section of the adjoint bundle. The boundary conditions imply that as $r=\vert{\bf x}\vert\rightarrow \infty$, $\vert\phi\vert\rightarrow 1-k/2r-\dots$ where $k$ is an integer, the magnetic charge. 
     
     The group $SO(3)$ of spatial rotations acts on the moduli space, rotating the complex structures of the hyperk\"ahler family, so rotations about a  fixed axis give us a circle action of the type we are considering. From a theorem of Donaldson we can identify the moduli space ${\mathcal M}_k$ as the space of rational functions $R_k$:
$$S(z)=\frac{p(z)}{q(z)}=\frac{a_0+a_1z+\dots+a_{k-1}z^{k-1}}{b_0+b_1z+\dots+z^k}.$$
The circle action is $S(z)\mapsto \lambda^{-2k}S(\lambda^{-1} z)$. From \cite{AH}  its twistor space is  of the form $Z=U\cup V$, $U\cong V\cong R_k\times \C$
where on $U\cap V=R_k\times \C^*$ we have the identification 
$$\tilde\zeta=\frac{1}{\zeta},\quad \tilde q\left(\frac{z}{\zeta^2}\right)=\frac{1}{\zeta^{2k}}q(z),\quad \tilde p\left(\frac{z}{\zeta^2}\right)=e^{-2z/\zeta}p(z)\,\, {\mathrm {mod}}\,\, q(z).$$
It is convenient to use local coordinates on the open set where the denominator has distinct roots $\beta_1,\dots,\beta_k$. Then since $p/q$ has degree $k$, $p(\beta_i)\ne 0$ and we use coordinates $\beta_i, p_i=p(\beta_i)$ on $U$ and $\tilde\beta_i,\tilde p_i$ on $V$ where
$$\tilde\beta_i=\beta_i/\zeta^2,\qquad \tilde p_i=\exp(-2\beta_i/\zeta)p_i.$$
The symplectic form along the fibres is $\sum_id\beta_i\wedge d\log p_i=\zeta^2\sum_id\tilde\beta_i\wedge d\log \tilde p_i$. Setting $\chi_i=\log p_i$ we have the same format as the Legendre transform, although there is no global $\C^k$-action here. If, with $q(z)=z^k+b_{k-1}z^{k-1}+\dots$,  we put
$$H(p(z)/q(z),\zeta)=\frac{2b_{k-2}-b_{k-1}^2}{\zeta}=-\frac{1}{\zeta}\sum_i\beta_i^2$$
then $\log \tilde p_i=\log p_i+\partial H/\partial \beta_i$. Moreover $H$ is homogeneous of degree one, so the holomorphic line bundle on the twistor space is defined by the transition function $\exp(-(\partial H/\partial\zeta)/2)$:
$$g_{UV}=\exp\left(\frac{2b_{k-2}-b_{k-1}^2}{2\zeta^2}\right).$$
When the monopole is centred $b_{k-1}=0$ and the resultant $R(p,q)=1$. Then we have simply $g_{UV}=\exp (b_{k-2}/\zeta^2)$.
\vskip .5cm
We can also give a differential-geometric description of the hyperholomorphic bundle, thanks to \cite{Hit4}, where the following formula is given for the K\"ahler potential of the moduli space of centred monopoles in the complex structure fixed by rotations about the unit vector ${\bf u}$:
\begin{equation}
\phi=\frac{4}
{(N+1)(N+2)}\frac{\vartheta^{(N+2)}(0)}{\vartheta^{(N)}(0)}-\frac{1}{3}Q({\bf u},{\bf u}).
\label{su2}
\end{equation}
Here $\vartheta$ is a theta function for the spectral curve (here assumed to be smooth) and $Q$ is the quadrupole moment in the asymptotic expansion of the Higgs field: 
$$\vert \phi \vert=1-\frac{k}{2r}-\frac{Q({\mathbf x},{\mathbf x})}{4r^5}+\dots$$
It is a trace zero symmetric rank 2 tensor. 

The spectral curve \cite{Hit1}, \cite{AH} is an algebraic curve whose equation is of the form $\eta^k+a_{k-1}(\zeta)\eta^{k-1}+\dots+a_0(\zeta)=0$ where $a_i(\zeta)$ is a polynomial of degree $2k-2i$, which is real in the sense that  $\zeta^{2k-2i}\overline{a_k(-1/\bar \zeta)}=a_k(\zeta)$. Then the quartic function  $2a_{k-2}(\zeta)-a_{k-1}^2(\zeta)=c_0+c_1\zeta+c_2\zeta^2-\bar c_1\zeta^3+\bar c_0 \zeta^4$ where $c_2$ is real. If ${\bf u}$ gives the complex structure $\zeta=0$ then $c_2=-Q({\bf u},{\bf u})/2$.  In the rational map description the numerator at complex structure $\zeta$ is given by setting $\eta=z$, so it is perhaps not surprising, given the transition function above, that  $Q$ appears in the formula.

Now consider the curvature $F=\omega_1+dd^c_1\mu$ of the hyperholomorphic line bundle. As we have seen, the K\"ahler form for $J$ is $\omega_2= -dd_2^c\mu$. But this is the complex structure invariant under rotation about an axis ${\bf v}$ orthogonal to ${\bf u}$. The theta function term in (\ref{su2}) is rotationally invariant, so we find 
$$F=-\frac{1}{3}dd_1^c(Q({\bf u},{\bf u})-Q({\bf v},{\bf v})).$$
Now under a rotation of an angle $\theta$, $c_2$ is replaced by 
\begin{equation}
c_2 +\frac{3}{2}\sin^2\theta (c_0+c_4-c_2)-\frac{3}{2}\sin \theta
\cos\theta (c_1-c_3)
\label{rot}
\end{equation}
so for the orthogonal complex structure defined by ${\bf v}$ we get
$$c_2 +\frac{3}{2} (c_0+\bar c_0-c_2)=-Q({\bf v},{\bf v})/2.$$
Hence 
$$Q({\bf u},{\bf u})-Q({\bf v},{\bf v})=-c_2+3(c_0+\bar c_0).$$
However $c_0= 2b_{k-2}-b_{k-1}^2$ and is holomorphic in the complex structure at $\zeta=0$, so we finally get
$$F=\frac{1}{6}dd_1^cQ({\bf u},{\bf u}).$$
 \subsection{Higgs bundle moduli spaces}
  Let $\Sigma$ be a compact Riemann surface and $P$ a principal $G$-bundle for $G$ a compact Lie group with complexification $G^c$. Given a  connection $A$ on $P$ and $\Phi$  a section of $\lie{g}\otimes K$ where $\lie g$ is the adjoint bundle, the Higgs bundle equations  are given by 
  \begin{equation}
  \bar\partial_A\Phi=0,\qquad F_A+[\Phi,\Phi^*]=0
  \label{Higgseq}
  \end{equation}
These are formally the equations for the zero set of a hyperk\"ahler moment map for the action of the gauge group, and as a consequence the moduli space ${\mathcal M}$ of solutions has a hyperk\"ahler structure \cite{Hit5}. There is also a natural circle action $(A,\Phi)\mapsto (A,e^{i\theta}\Phi)$ which fixes one K\"ahler form and rotates the other two. 

The different complex structures, parametrized by  $\zeta\in \CP^1$, can be viewed by considering the operators $D_+:\Omega^0(\lie{g})\rightarrow \Omega^{0,1}(\lie{g})$ and $D_-:\Omega^0(\lie{g})\rightarrow \Omega^{1,0}(\lie{g})$ defined by 
$$D_+=\bar\partial_A+\zeta\Phi^*,\qquad D_-=\partial_A+\tilde\zeta\Phi$$
For $\zeta\ne 0,\infty$ and $\tilde \zeta=\zeta^{-1}$, then $D_+=\nabla^{0,1},D_-=\nabla^{1,0}$ for a flat $G^c$-connection $\nabla$. When $\zeta=0$, $D_+$ is the $\bar\partial$-operator $\bar\partial_A$ which defines a holomorphic structure on the principal $G^c$-bundle and for which the equations (\ref{Higgseq}) say that $\Phi$ is holomorphic and satisfies a stability condition. The gauge-theoretical moduli space thus has different interpretations, either as a moduli space of Higgs pairs $(\bar\partial_A,\Phi)$ when $\zeta=0$ or $\infty$ or as a moduli space of flat connections otherwise. 

The twistor space $Z$ is the union of two open sets $U,V$ given by $\zeta\ne \infty$ and $\zeta\ne 0$ respectively. Each one is actually algebraic, a moduli space of  {\it $\lambda$-connections} in Deligne-Simpson's terminology \cite{Simp}: for the general linear group these are  holomorphic differential operators $D: E\rightarrow \Omega^{1,0}(E)$ on a vector bundle $E$ such that $D(fs)=\lambda \partial f\otimes s+fDs$. On $U$, $D_+=\bar\partial_A+\lambda\Phi^*$ defines the holomorphic structure on $E$ and the $\lambda$-connection is $D=\lambda\partial_A+\Phi$. For a fixed holomorphic structure these form a vector space whose projective space has a distinguished hyperplane $\lambda=0$. The affine space which is the complement is the space of holomorphic  connections on $E$.

An  elliptic operator $D$ on a compact manifold  has finite dimensional kernel and cokernel and one can associate to this a  determinant line 
$$L(D)=(\Lambda^{top}\!\ker D)^*\otimes \Lambda^{top}\!\coker D.$$
and in particular for a $\bar\partial$-operator on sections of a vector bundle on a Riemann surface. It behaves well in families even though the dimensions may jump (see \cite{freed}, \cite{Q}) and defines a  line bundle $L$. In particular, the $\bar\partial$-operator for a $\lambda$-connection ($D_+$ with $\lambda=\zeta$)  defines  a determinant line bundle $L_+$ over $U$. There is a similar line bundle $L_-$ on $V$ defined by $D_-$. 

For the  moduli space of Higgs bundles $H^2({\mathcal M},\Z)$ has a single generator, so all determinant bundles are powers of one. To deal with a general group it is convenient to take the bundle given by the adjoint representation and to tensor with a line bundle $K^{1/2}$, whose square is the canonical bundle, a so-called theta-characteristic. We shall see this later, but it is immaterial for the following:  

\begin{prp} \label{Lplusminus}The holomorphic line bundle $L_Z$ on the twistor space $Z$  defining the hyperholomorphic line bundle of ${\mathcal M}$ is isomorphic to $L_+$ on $U$ and $L^*_-$ on $V$.
\end{prp}
\begin{prf} On the infinite-dimensional affine space ${\mathcal A}$ of $\bar\partial$-operators on the  bundle ${\lie g}\otimes K^{1/2}$, the determinant line bundle $L$ has a natural connection defined by the Quillen metric \cite{Q}. The curvature of the Quillen metric is a multiple  of the Hermitian form $\int_{\Sigma}\tr da da^*$. On the flat hyperk\"ahler manifold ${\mathcal A}\times H^0(\Sigma, \lie{g}\otimes K)$ the hyperholomorphic bundle has curvature 
$$F=\frac{1}{2i} \int_{\Sigma}(\tr da da^*-\tr d\Phi d\Phi^*).$$
Recall from (\ref{herm}) that in a finite-dimensional flat space with Hermitian  form  $\langle z,\bar z\rangle=\sum_iz_i\bar z_i$ if $\log h_U=\langle z,\bar z\rangle-\langle w,\bar w\rangle+\zeta\langle\bar z,\bar w\rangle+\bar\zeta \langle z, w \rangle $ then $\bar\partial_Z \partial_Z \log h_U=\langle dz, d\bar z\rangle-\langle dw,d\bar w \rangle=2iF$. Choosing a basepoint $D_0\in {\mathcal A}$ we can do the same thing on ${\mathcal A}\times H^0(\Sigma, \lie{g}\otimes K)$ as the circle acts trivially on ${\mathcal A}$.

Now note that 
$$\langle z+\zeta \bar w, \bar z+\bar\zeta w\rangle-(1+\zeta\bar\zeta)\langle w,\bar w\rangle=\langle z,\bar z\rangle-\langle w,\bar w\rangle+\zeta\langle\bar z,\bar w\rangle+\bar\zeta \langle z, w \rangle.$$
Since $z_i+\zeta \bar w_i$ are, for $\zeta\ne \infty$, holomorphic functions on $Z$ we have 
$$2iF= \bar\partial_Z \partial_Z \log h_U=\langle d(z+\zeta \bar w)\wedge d(\bar z+\bar \zeta w)\rangle- \partial_Z \bar\partial_Z((1+\zeta\bar\zeta)\langle w,\bar w\rangle). $$

In the present  context this means that taking the  Quillen metric for the determinant bundle of $D_+=\bar\partial_A+\zeta\Phi^*$ and rescaling it by $\exp [(1+\zeta\bar\zeta)\int_{\Sigma}\tr\Phi\Phi^*]$ we have a metric whose curvature is $F$, and hence the line bundle $L_Z$ on $U$ is isomorphic to the determinant line bundle $L_+$. A similar argument gives $L^*_-$ on the open set $V$.
\end{prf}

To describe the line bundle $L_Z$ we now need not a transition function, but instead an isomorphism between $L_+$ and $L_-^*$ on $U\cap V$.

Note that $\Phi=0$ is a component of the fixed point set of the circle action, and this is the moduli space ${\mathcal N}$ of stable  bundles, or flat connections with holonomy in the compact group $G$. The tangent space at  a point in ${\mathcal N}$ is naturally isomorphic to $H^1(\Sigma,\lie{g})$. By Serre duality the cotangent space is $H^0(\Sigma, \lie{g}\otimes K)$, where the Higgs field lies, and so there is a neighbourhood of ${\mathcal N}\subset \mathcal{M}$ which is identified with the cotangent bundle of ${\mathcal N}$, with the circle action given by scalar multiplication in the fibres. The hyperk\"ahler metric is therefore, using Feix's result,  the unique hyperk\"ahler extension of the K\"ahler metric on the moduli space of stable bundles. As in Section \ref{cotangent}, we can in theory describe the line bundle $L_Z$ by transition functions involving the K\"ahler potential of this metric. 

The K\"ahler form on ${\mathcal N}$ is known to be a multiple of the curvature of the Quillen metric.  On the infinite-dimensional  space ${\mathcal A}$ the K\"ahler potential is, as we saw above,  just the Hermitian  form, but because of the choice of base-point this is not gauge-invariant. Quillen gave an alternative description involving the zeta-function regularized determinant 
$$\mathrm{det}_{\zeta}(\bar\partial_A^*\bar\partial_A)=\exp (-\zeta'(0))$$ 
of the composition of an operator $\bar\partial_A$ with its adjoint. Here $\zeta(s)$ is the analytic continuation of $\zeta(s)=\sum_{\lambda\ne 0}\lambda^{-s}$ for the eigenvalues $\lambda$ of $\bar\partial_A^*\bar\partial_A$. This depends on a choice of metric, but only up to a factor which depends on the metric alone. 

The operator $\bar\partial_A:\Omega^0(\lie{g}\otimes K^{1/2})\rightarrow \Omega^{0,1}(\lie{g}\otimes K^{1/2})$ has the property that  $\dim \ker \bar\partial_A=\dim\coker\bar\partial_A$.  It then follows  (see \cite{freed}) that  there is a canonical determinant section $\sigma(\bar\partial_A)$ of the determinant line bundle $L$. (In fact  if the theta characteristic is odd and the dimension of $G$ is odd this always vanishes, but there are choices for which, as $A$ varies, it is generically non-zero.) The Quillen metric is then defined by the property that 
\begin{equation}
\Vert \sigma(\bar\partial_A)\Vert^2=\mathrm{det}_{\zeta}(\bar\partial_A^*\bar\partial_A).
\label{quill}
\end{equation}
The determinant section $\sigma(\bar\partial_A)$ vanishes when $\bar\partial_A$ has a non-zero kernel, but the behaviour of the zeta-function and the definition of $\sigma(\bar\partial_A)$ as in   \cite{freed}, show that the metric, a section of $L^*\otimes \bar L^*$ is everywhere non-vanishing. Since the curvature form for this bundle is the K\"ahler form on ${\mathcal N}$, the local trivialization of $L$ defined by $\sigma(\bar\partial_A)$ gives from (\ref{quill}) a K\"ahler potential $f/2=\log\mathrm{det}_{\zeta}(\bar\partial_A^*\bar\partial_A)$. 
 To implement the construction of Section \ref{cotangent} we need to understand the complexification of $\exp f/2=\mathrm{det}_{\zeta}(\bar\partial_A^*\bar\partial_A)$ and then identify it as an isomorphism between the determinant line bundles.
 
Recall that  the construction of Section \ref{cotangent}  uses the diffeomorphism $\phi_+:N^c\times \C^*\rightarrow V_+^*$. In the Higgs bundle situation  the K\"ahler manifold $N$ is the moduli space of flat $G$-connections ${\mathcal N}$ and its complexification $N^c$ we view as  the moduli space of flat $G^c$-connections and $V_+^*$  the moduli space of $\lambda$-connections $D$. For fixed $\lambda\in \C^*$, if we put $\nabla^{0,1}=\bar\partial$ and $\nabla^{1,0}=\lambda^{-1}D$ then $\nabla=\nabla^{0,1}+\nabla^{1,0}$ is  a flat $G^c$-connection and this defines $\phi_+^{-1}$. 

Given the flat $G^c$-connection $\nabla$ we can define two operators
$$ D_+:\Omega^0(\lie{g}\otimes K^{1/2})\stackrel{\nabla^{0,1}}\rightarrow \Omega^{0,1}(\lie{g}\otimes K^{1/2})\quad D_-: \Omega^0(\lie{g}\otimes \bar K^{1/2})\stackrel{\nabla^{1,0}}\rightarrow \Omega^{1,0}(\lie{g}\otimes \bar K^{1/2}).$$
These have determinant bundles $L^+,L^-$ as considered in Proposition \ref{Lplusminus}, and determinant sections $\sigma_+,\sigma_-$. A metric in $\Sigma$ defines a section $h$ of $K\bar K$ and then $\Delta = h^{-1/2}D_- h^{-1/2}D_+$ is an operator which specializes to $\bar\partial_A^*\bar\partial$ in the case of a $G$-connection. 
Moreoever it has the same principal symbol and so satisfies the Agmon-Nirenberg condition on its spectrum which guarantees a well-defined (complex) zeta-function determinant. The definition of $\sigma_+,\sigma_-$ and the properties of the zeta-function show that $\sigma_+\sigma_-/\det_\zeta \Delta$ is a non-vanishing section of $L_+\otimes L_-$ and this defines the required identification $L_+\cong L_-^*$. In fact we know from Proposition \ref{Lplusminus} that such an isomorphism must exist and we also know that it is uniquely determined by the analytic continuation of the K\"ahler potential, but this combination of determinants gives it a more concrete form. 

\begin{rmk} Two recent papers \cite{Pio}, \cite{Neitz}  discuss the hyperholomorphic bundle in terms of wall-crossing, the first using the existence of a circle action, the second apparently not. It nevertheless uses Higgs bundles, but with singularities, and these  generally do not have such an action. There are however determinant lines, so it seems  possible that the above description is concerned with  the same hyperholomorphic bundle. 
\end{rmk}
We do not know an explicit form for the Higgs bundle metrics for genus $g>1$ but we can observe the description above in the case of an elliptic curve, which we do next.

\subsection{Higgs bundles on an elliptic curve}
We shall consider here the moduli space ${\mathcal M}$ of $U(1)$-Higgs bundles on an elliptic curve $\Sigma$ with modulus $\tau=x+iy, y>0$. The moduli space of flat $U(1)$-bundles is a flat torus, and then ${\mathcal M}$ is just $T^2\times \R^2$ with a flat metric, but we shall look at it from the point of view of the previous section and then compare with what we know of the hyperholomorphic bundle in the flat case. Since the adjoint bundle is trivial in this case we consider the determinant line for the vector representation, i.e. we look first at a flat $U(1)$-line bundle ${\mathcal L}$ and the determinant line of the operator
$\bar\partial_A:\Omega^0({\mathcal L})\rightarrow \Omega^{0,1}({\mathcal L}).$

The zeta-function determinant was calculated in \cite{RS}: for a non-trivial character given by $A\mapsto e^{2\pi i a}, B\mapsto e^{2\pi i b}$ for $0\le a,b <1$
$$-\zeta'(0)=\left\vert e^{\pi i a^2 \tau}\frac{\vartheta(b-\tau a,\tau)}{\eta(\tau)}\right\vert^2=e^{\pi i a^2 (\tau-\bar\tau)}\frac{\vartheta(b-\tau a,\tau)\overline{\vartheta(b-\tau a,\tau)}}{\vert\eta(\tau)\vert^2}
$$
The theta function $\vartheta$ here has the expansion
$$\vartheta(z,\tau)=-2q^{1/4}\sin (\pi z) \prod_{m=1}^{\infty}(1-q^{2m})(1-2\cos (2\pi z) q^{2m}+q^{4m})  $$
where $q=\exp \pi i \tau$.
It is invariantly to be thought of as a  section  of the determinant line bundle -- it vanishes where $a=b=0$, where the holomorphic line bundle ${\mathcal L}$ is trivial and  has a non-zero section,  and is thus a multiple of the determinant section $\sigma(\bar\partial_A)$. It is only a function as expressed here when lifted to the universal covering of the torus. In differential-geometric notation, if $du$ is the non-vanishing 1-form on $\Sigma$ with periods $1,\tau$ and the flat connection is $d+\alpha d\bar u-\bar \alpha du$ for $\alpha\in \C$ then $2\pi i(b-\tau a)=(\tau-\bar \tau)\alpha=2iy\alpha$. 

When we  complexify  $\zeta'(0)$, the variable $b-\tau a$ and its conjugate become independent variables $z,\tilde z$ and then the expansion shows that
\begin{equation}
-\zeta'(0)=\exp(\pi  (z-\tilde z)^2/ 2y)\frac{\vartheta(z,\tau){\vartheta(\tilde z,-\bar\tau)}}{\vert\eta(\tau)\vert^2}
\label{complextheta}
\end{equation}

A solution to the Higgs bundle equations in the abelian case is given by a flat unitary connection $d+\alpha d\bar u-\bar \alpha du$ and a constant Higgs field $\Phi=\beta du$. Together they define the flat $\C^*$-connection $\nabla_A+\zeta\Phi^*+\zeta^{-1}\Phi$ which is 
$$\nabla=d+\alpha d\bar u-\bar \alpha  du+\zeta \bar \beta d\bar u+\zeta^{-1}\beta du$$
so for the $\nabla^{0,1}$-operator of this we get $2\pi i z=2\pi i(b-\tau a)=2iy(\alpha+\zeta\bar \beta).$ For the conjugate structure $2\pi i \tilde z=2iy(\bar \alpha-\zeta
^{-1}\beta)$. 

Consider now the interpretation of  (\ref{complextheta}) on the universal covering of the twistor space, which is just the flat twistor space $\C^2(1)$. A Higgs bundle defines a real section, and using holomorphic coordinates $(v,\xi,\zeta)$ as in Example  \ref{flat}, this section is $v=(\alpha+\zeta\bar\beta)y, \xi =(\beta-\zeta\bar\alpha)y$.  Then with $z=v/\pi, \tilde z=-\xi/\zeta\pi$  (\ref{complextheta}) defines a local transition function for a  holomorphic line bundle. The two theta function factors and the constant  eta-function term are local trivializations of the determinant bundles $L_+,L_-$ and so the isomorphism between them is given by 
$$\exp(\pi  (z-\tilde z)^2/ 2y)=\exp( (v+\xi/\zeta)^2/ 2\pi y).$$ 
As in Section \ref{Leg} $\exp v^2/2\pi y$ and $\exp \xi^2/\zeta^22\pi y$ are changes of local trivialization on the open sets $U,V$ and with respect to these the isomorphism is defined by 
$$\exp(v\xi/ \pi y \zeta).$$
 This is essentially the usual flat space description of the line bundle $L_Z$, the factor $\pi y$ appearing in order to   implement the integrality condition on the cohomology class.

\section{The Quaternionic K\"ahler manifold}
\subsection{Quaternionic K\"ahler geometry}
A quaternionic K\"ahler manifold is a Riemannian manifold $M^{4k}$ whose holonomy is contained in $Sp(k)\cdot Sp(1)$. A hyperk\"ahler manifold has holonomy in $Sp(k)\subset Sp(k)\cdot Sp(1)$ but we usually distinguish between the two notions: the scalar curvature is a non-zero constant for quaternionic K\"ahler. The $\pm Sp(1)$-factor defines a principal $SO(3)$-bundle $S$. This is the frame bundle for a bundle of quaternion algebras which act on the tangent bundle of $M$, or equivalently a rank $3$ bundle of $2$-forms $\omega_1,\omega_2,\omega_3$ which satisfy the algebraic relations of hyperk\"ahler forms.  

There is a link between the two geometries. The connection on $S$ and an orthonormal basis for $\lie{so}(3)$ define three $1$-forms $\theta_1,\theta_2,\theta_3$ on $S$. The components $K_{23}=d\theta_1-\theta_2\wedge\theta_3$ etc. of the curvature of the connection are of the form $K_{23}=c\omega_1$ where $c$ is essentially the scalar curvature. On the $4k+4$-manifold $S\times \R^+$  we have three closed $2$-forms $\varphi_i=d(t\theta_i)$. Now $T(S\times \R^+)=H\oplus V$ where $H$ is the horizontal subbundle defined by $\theta_i=0$ and $dt=0$ and restricted to  this subspace we have  $\varphi_i=tc\omega_i$. On $V$ we have $\varphi_1=dt\wedge\theta_1+t^2\theta_2\wedge\theta_3$ etc. so together these define an $Sp(k+1)$-structure on the tangent space if $c>0$ and $Sp(1,k)$ if $c<0$. This is in fact hyperk\"ahler and then $S\times \R^+$ is called the {\it Swann bundle} \cite{swann} or hyperk\"ahler cone of the quaternionic K\"ahler manifold $M$. 

The quotient space $S/SO(2)$ is the unit sphere bundle in the bundle of imaginary quaternions and so is a bundle of complex structures on the tangent bundle of $M$. Using the connection, we have a splitting  $T(S/SO(2))=H\oplus V$ and, as in the hyperk\"ahler case, we can introduce an almost complex structure $(I_{\mathbf x},I)$. This is integrable and defines the twistor space of $M$.

We can see this via the Swann bundle: the $SO(3)$ action on $S$ rotates the complex structures $I,J,K$ and so $SO(2)\subset SO(3)$ fixes $I$, say. Then $\C^*=SO(2)\times \R^+$ acts $I$-holomorphically on $S\times \R^+$ and the twistor space is the quotient. It is a complex manifold of dimension $2k+1$ but has an extra structure: if $X$ is the holomorphic vector field generated by the $\C^*$-action on $S\times \R^+$, let $X^{\perp}$ be the symplectic-orthogonal subbundle of the tangent bundle with respect to $\omega_2+i\omega_3$. Then $X^{\perp}/X$  descends to the quotient as a rank $2k$ distribution which is a contact structure. 

Recall that if $L$ is the line bundle on the twistor space which is the quotient by this codimension one distribution then an equivalent description is via a section $\alpha$ of $T^*(L)=T^*\otimes L$. For this to be a contact form we require $\alpha\wedge(d\alpha)^k$ to be everywhere non-vanishing which implies $L^{k+1}\cong K^*$, which we write $L=K^{-1/(k+1)}$.

\begin{ex} If $M=\KP^k$ then its twistor space is $\CP(\C^{2k+2})$ where $\C^{2k+2}$ has a non-degenerate skew form. Then ${\mathcal O}(-1)^{\perp}/{\mathcal O}(-1)$ defines the contact structure $\alpha$ which is a section of $T^*(2)$.
\end{ex}

As a complex manifold, the twistor space $Z$ has the following features \cite{Sal}:
\begin{itemize}
\item
a complex $2k+1$-manifold
\item
a holomorphic section  $\alpha$ of $T^*\otimes K^{-1/(k+1)}$ such that $\alpha\wedge(d\alpha)^k\ne 0$
\item
a family of  rational curves with normal bundle $\C^{2k}(1)$ and on which $\alpha$ is nonzero 
\item
a real structure preserving this data.
\end{itemize}
With this information we can reconstruct the quaternionic K\"ahler manifold $M$.

\subsection{The correspondence: twistor approach}
Let us return now to the result of Section \ref{twist}. On the twistor space $Z$ of a hyperk\"ahler manifold with circle action we constructed  a principal $\C^*$-bundle $P$ and lifted the action. It is convenient to assume that the holomorphic vector field defined by the circle generates a $\C^*$-action and then, removing the fixed points,  we let  $\hat Z$ be the quotient.  In fact one only needs a local holomorphic extension of the circle action to achieve this.  
We shall show that $\hat Z$ has a contact structure $\alpha$ invariant by  the induced $\C^*$-action. 

Let $\hat p:P\rightarrow \hat Z$ be the quotient map by the action. Then since the tangent bundle along the fibres is trivial $K_P\cong \hat p^*K_{\hat Z}$. By the same token, $K_P\cong p^*K_Z\cong {\mathcal O}(-2(k+1))$ and 
\begin{equation}
\hat p^*K_{\hat Z}^{-1/(k+1)}\cong p^*{\mathcal O}(2). 
\label{isobundle}
\end{equation}
We shall continue the notation ${\mathcal O}(2)$ since it is of degree $2$ on each twistor line, but it is canonically determined on $\hat Z$. A contact form $\alpha$ is then a  section of $T^*_{\hat Z}(2)$.

So if $\alpha$ is a contact form on $\hat Z$, $\hat p^*\alpha$ should be  a section of $T^*_P(2)$, invariant by the $\C^*$-action generated by the vector field $\tilde Y$, and such that $i_{\tilde Y}\hat p^*\alpha=0$. This is what we shall define.

In Proposition \ref{mero} we defined a meromorphic connection on $L_Z$, so there is a meromorphic connection 1-form $A$ on $P$,  with simple poles at $\zeta=0,\infty$. This defines a {\it holomorphic} section $\zeta A$ of $T^*_P(2)$. It is invariant under the circle action, so to define $\alpha$ on $\hat Z$ we need to check that $i_{\tilde Y}A=0$. Locally on $P$ the connection form is 
$$A=dt-\frac{1}{i\zeta} (i_{Y_U}\omega_U-f_Ud\zeta)$$
and the lifted vector field $\tilde Y=Y-f_U V$ where $V$ is the vertical vector field given by the $\C^*$-principal bundle action. So
$i_{\tilde Y}A=-f_U-i_{Y_U}i_{Y_U}\omega_U/i\zeta+f_U=0$, and $\zeta A=\hat p^*\alpha$ for a section of $T^*_{\hat Z}(2)$.

We need to show that this has the nondegeneracy property of a contact structure. But 
$\zeta A \wedge d(\zeta A)^k=\zeta^{k+1}A\wedge {\mathcal F}^k$ where ${\mathcal F}$ is the curvature of the meromorphic connection, and we showed in Proposition \ref{mero} that this was nondegenerate restricted to the generic fibres of $p:Z\rightarrow \CP^1$. Since the connection form $A$ is non-zero on the fibres of $P\rightarrow Z$ we see that $\alpha\wedge (d\alpha)^k$ is non-zero for $\zeta\ne 0,\infty$. But it is a section of the trivial bundle on a  twistor space and is therefore constant hence non-zero everywhere.

We have thus shown that $\hat Z$ is a holomorphic contact manifold. The bundle $P$ is trivial on the twistor lines of $Z$ so they lift to $P$ and project to $\hat Z$ to define the twistor lines there.

\subsection{The converse}
The reverse construction, from a quaternionic K\"ahler manifold with $S^1$-action to a hyperk\"ahler manifold can be done in two ways. For the first, we observe how to invert the above twistorial construction.  

\noindent 1. The section $s$ of ${\mathcal O}(2)$ on $Z$ which vanishes on the divisor $D_0+ D_{\infty}$ is invariant by the circle action and so we have a section $\hat s$ on $\hat Z$ vanishing on the divisor $\hat D_0+\hat D_{\infty}$. Now consider $\zeta$ as a meromorphic function on $Z$, with a zero on $D_0$ and a pole on $D_{\infty}$. Under the $\C^*$-action it transforms as $\zeta\mapsto\lambda\zeta$. Pulled back to $P$, this becomes a meromorphic section of the line bundle $\hat L$ over $\hat Z$ defined by $P$. Put another way, $\hat D_0-\hat D_{\infty}$ is  the divisor class of the line bundle $\hat L$. This provides the key to the correspondence.

So suppose $\hat M$ is a quaternionic K\"ahler manifold with a circle action. The action generates a holomorphic vector field $W$ on its twistor space $\hat Z$ and then $i_W\alpha$ is a section  $\mu$ of ${\mathcal O}(2)$. This is the twistorial version of the quaternionic K\"ahler moment map of \cite{Gal}. The section vanishes at two points on a generic twistor line and we assume that the divisor consists of two components $\hat D_0$ and $\hat D_{\infty}$, interchanged by the real structure. Let $L$ be the line bundle, of degree zero on each twistor line, defined by the divisor $\hat D_0-\hat D_{\infty}$, and let $\hat p:P\rightarrow \hat Z$ be the corresponding principal $\C^*$-bundle. Then on $P$, $\hat p^*L$ is trivial and $p^*(\hat D_0-\hat D_{\infty})$ is the divisor of a meromorphic function $\zeta$, in other words we have a holomorphic map $P\rightarrow \CP^1$ with $\zeta=0,\infty$ being the divisors $\hat p^*\hat D_0,\hat p^*\hat D_{\infty}$. There is a lift of the action, giving a vector field $\tilde W$ on $P$, such that  the function $\zeta$ is invariant and then the map descends to the quotient $Z$, and describes the fibration $p:Z\rightarrow \CP^1$. 

We need the section $\omega$ of $\Lambda^2T^*_F(2)$ on $Z$ and this is obtained from differentiating $\hat p^*\alpha\in T^*_P(2)$ in the fibre directions.  Since $i_{W}\alpha=\zeta$ and ${\mathcal L}_{W}\alpha=0$, we have $i_{\tilde W}d \hat p^*\alpha =-d\zeta$ and so along a fibre $\zeta= const.$, $d \hat p^*\alpha$ descends to the quotient $Z$ as the required section $\omega$. 
\vskip .25cm
\noindent 2. The second procedure, which can be implemented in both the differential-geometric and twistorial approaches, is to use the Swann bundle. The circle action on the quaternionic K\"ahler manifold $\hat M$ has a canonical lift to the hyperk\"ahler Swann bundle and then the corresponding hyperk\"ahler manifold $M$ is  the hyperk\"ahler quotient. 

We proceed as follows, describing first the twistor space of the Swann bundle (see \cite{Sal1} for the original description). Let $q:Q\rightarrow \hat Z$ denote the principal $\C^*$-bundle for the line bundle ${\mathcal O}(2)$ and consider the associated bundle  $Q\times_{\C^*}\C^2\backslash\{0\}/\{\pm 1\}$ where the action of $\lambda\in \C^*$ on $\C^2\backslash\{0\}/\{\pm 1\}$ is induced from $\lambda\cdot(z_0,z_1)\mapsto (\lambda^{1/2} z_0,\lambda^{1/2} z_1)$. Then the homogeneous coordinates $(z_0,z_1)$ define a projection to $\CP^1$. Moreover $q^*\alpha$ is a 1-form on $Q$ whose derivative is a symplectic form (the canonical symplectic form on the bundle of contact elements). Adding $dz_0\wedge dz_1$ and taking the quotient by $\C^*$ gives a section of $\Lambda^2T^*_F(2)$, making it the twistor space of the hyperk\"ahler structure on the Swann bundle.  

Now take the lifted action, which has a moment section of ${\mathcal O}(2)$ defined by $\mu+z_0z_1$. The twistor space for the hyperk\"ahler quotient is 
the zero set of this modulo the $\C^*$-action. On this subvariety we have the relation $\mu=-z_0z_1$ and so $\mu$ vanishes on two divisors given by $z_0=0$ and $z_1=0$. Since $z_0/z_1$ is a meromorphic function the two divisors are linearly equivalent. The action $(z_0,z_1)\mapsto (\nu z_0,\nu^{-1}z_1)$ expresses $\mu+z_0z_1=0$ as a principal $\C^*$-bundle over $\hat Z$ and so we have recovered the first description. 

\begin{rmk} The two divisors $\hat D_0,\hat D_{\infty}$ in the twistor space $\hat Z$ of the quaternionic K\"ahler manifold define sections of $\hat Z\rightarrow \hat M$, interchanged by the real structure. These can be interpreted as complex structures $\pm I$ on $\hat M$. These are examples of manifolds with a torsion-free connection with holonomy in $SL(k,\K)\cdot U(1)$. When $k=1$ this group is $U(2)$ and the structure is that of a scalar-flat K\"ahler metric but in higher dimensions it is a non-metric geometry. This seems to have been little studied apart from \cite{Joy} where they are called {\it quaternionic complex} manifolds. Whether this structure or the quaternionic K\"ahler one is more important globally in the correspondence studied here is a question we leave till another time. 
\end{rmk}

\subsection{An example}
We shall consider only a single example of a hyperk\"ahler manifold with a circle action -- the case where $M$ is the Eguchi-Hanson metric on $T^*S^2$. This is the cotangent bundle case with the natural action on the fibres. The induced metric on the zero section is $SO(3)$-invariant and hence a multiple of the standard  metric on $S^2$: the metric is thus a case of the hyperk\"ahler extension of Feix. On the other hand, the metric appears also as a hyperk\"ahler quotient of flat space by a circle action, and the twistor space  has a description as the quotient of an open set in 
 \begin{equation}
\{(v,\xi)\in V(1)\oplus V^*(1)\rightarrow \CP^1:\langle v,\xi\rangle =\zeta\}.
\label{Zspace}
\end{equation}
where $V$ is a $2$-dimensional vector space, by the action   $(v,\xi)\mapsto (\lambda v,\lambda^{-1}\xi)$. Since this commutes with the $U(2)$-action on $V$, the corresponding hyperholomorphic line bundle must be $U(2)$-invariant. There is up to a multiple only one invariant 2-form on $S^2$, so by the uniqueness of the hyperholomorphic extension this line bundle obtained through a quotient is also the one constructed here. Thus (\ref{Zspace}) with the $\C^*$-action is the holomorphic principal bundle $P$ for Eguchi-Hanson.

 Consider the lift  $(v,\xi,\zeta)\mapsto(\nu v,\nu\xi, \nu^2\zeta)$ of the geometrical action on $M$. Then the equation $\langle v,\xi\rangle =\zeta$ determines $\zeta$ in terms of $v\in V$ and $\xi\in V^*$. The quotient is therefore an open set in $\CP^3=\CP(V\oplus V^*)$ which is the twistor space of $S^4$. The quaternionic K\"ahler metric on $\hat M$ is then the  standard metric on the sphere. 

Strictly speaking, the quaternionic K\"ahler transform is the complement of a circle in $S^4$ because we have to remove fixed points of the circle action on the hyperk\"ahler side. We can also see this on the quaternionic K\"ahler side: the moment section of ${\mathcal O}(2)$ in $\CP(V\oplus V^*)$ is defined by  $\langle v,\xi\rangle$. This vanishes on a quadric,  isomorphic to $\CP^1\times \CP^1$ which is connected, but removing $\RP^1\subset \CP^1$ in the second factor produces two components -- our divisors $\hat D_0,\hat D_{\infty}$ -- and each of these has the conformally flat scalar-flat K\"ahler metric which is  the product of constant curvature metrics on $H^2\times S^2$

There is a choice involved in lifting the geometric  action on the twistor space -- any two lifts differ by a power of the principal bundle action.  A different choice of lifting is given by  $(v,\xi,\zeta)\mapsto(\nu^{n+1} v,\nu^{1-n}\xi, \nu^2\zeta)$. The quaternion K\"ahler manifolds produced here were considered by Haydys \cite{Hay} who in turn relates them to quotient constructions of Galicki and Lawson \cite{GL}.

\vskip 1cm
 Mathematical Institute, 24-29 St Giles, Oxford OX1 3LB, UK
 
 hitchin@maths.ox.ac.uk

 \end{document}